\documentclass[10pt]{article}
\usepackage{tikz}
\usetikzlibrary{calc}
\usepackage{amsmath, upgreek}
\usepackage{amssymb}
\usepackage{color}
\usepackage{amscd}
\usepackage{xspace}
\usepackage{verbatim}
\usepackage{graphicx}
\renewcommand{\theequation}{\thesection.\arabic{equation}
}

 \usepackage{cite}

\begin{document}


\newcommand{\txt}[1]{\;\text{ #1 }\;}
\newcommand{\tbf}{\textbf}
\newcommand{\tit}{\textit}
\newcommand{\tsc}{\textsc}
\newcommand{\trm}{\textrm}
\newcommand{\mbf}{\mathbf}
\newcommand{\mrm}{\mathrm}
\newcommand{\bsym}{\boldsymbol}
\newcommand{\scs}{\scriptstyle}
\newcommand{\sss}{\scriptscriptstyle}
\newcommand{\txts}{\textstyle}
\newcommand{\dsps}{\displaystyle}
\newcommand{\fnz}{\footnotesize}
\newcommand{\scz}{\scriptsize}
\newcommand{\be}{\begin{equation}}
\newcommand{\bel}[1]{\begin{equation}\label{#1}}
\newcommand{\ee}{\end{equation}}
\newcommand{\eqnl}[2]{\begin{equation}\label{#1}{#2}\end{equation}}
\newcommand{\barr}{\begin{eqnarray}}
\newcommand{\earr}{\end{eqnarray}}
\newcommand{\bars}{\begin{align*}}
\newcommand{\ears}{\end{align*}}
\newcommand{\nnu}{\nonumber \\}
\newtheorem{subn}{\name}
\renewcommand{\thesubn}{}
\newcommand{\bsn}[1]{\def\name{#1}\begin{subn}}
\newcommand{\esn}{\end{subn}}
\newtheorem{sub}{\name}[section]
\newcommand{\dn}[1]{\def\name{#1}}   
\newcommand{\bs}{\begin{sub}}
\newcommand{\es}{\end{sub}}
\newcommand{\bsl}[1]{\begin{sub}\label{#1}}
\newcommand{\bth}[1]{\def\name{Theorem}
\begin{sub}\label{t:#1}}
\newcommand{\blemma}[1]{\def\name{Lemma}
\begin{sub}\label{l:#1}}
\newcommand{\bcor}[1]{\def\name{Corollary}
\begin{sub}\label{c:#1}}
\newcommand{\bdef}[1]{\def\name{Definition}
\begin{sub}\label{d:#1}}
\newcommand{\bprop}[1]{\def\name{Proposition}
\begin{sub}\label{p:#1}}

\newcommand{\aand}{\quad\mbox{and}\quad}
\newcommand{\M}{{\cal M}}
\newcommand{\A}{{\cal A}}
\newcommand{\B}{{\cal B}}
\newcommand{\I}{{\cal I}}
\newcommand{\J}{{\cal J}}
\newcommand{\D}{\displaystyle}
\newcommand{\RR}{ I\!\!R}
\newcommand{\C}{\mathbb{C}}
\newcommand{\R}{\mathbb{R}}
\newcommand{\Z}{\mathbb{Z}}
\newcommand{\N}{\mathbb{N}}
\newcommand{\T}{{\rm T}^n}
\newcommand{\cuad}{{\sqcap\kern-.68em\sqcup}}
\newcommand{\abs}[1]{\mid #1 \mid}
\newcommand{\norm}[1]{\|#1\|}
\newcommand{\equ}[1]{(\ref{#1})}
\newcommand\rn{\mathbb{R}^N}
\newcommand\s{s}
\renewcommand{\theequation}{\thesection.\arabic{equation}}
\newtheorem{definition}{Definition}[section]
\newtheorem{theorem}{Theorem}[section]
\newtheorem{proposition}{Proposition}[section]
\newtheorem{example}{Example}[section]
\newtheorem{proof}{proof}[section]
\newtheorem{lemma}{Lemma}[section]
\newtheorem{corollary}{Corollary}[section]
\newtheorem{remark}{Remark}[section]
\newcommand{\bremark}{\begin{remark} \em}
\newcommand{\eremark}{\end{remark} }
\newtheorem{claim}{Claim}

\newcommand{\cH}{\mathcal{H}}
\newcommand{\cL}{\mathcal{L}}
\newcommand{\cO}{\mathcal{O}}
\newcommand{\cA}{\mathcal{A}}
\newcommand{\cJ}{\mathcal{J}}
\newcommand{\cN}{\mathcal{N}}
\newcommand{\cF}{\mathcal{F}}
\newcommand{\cE}{\mathcal{E}}
\newcommand{\cK}{\mathcal{K}}
\newcommand{\BBE}{\mathbb{E}}

 \begin{center} {\large \bf  Counting function and lower bound for Dirichlet eigenvalues \\[2mm]
   of the $m$-order logarithmic Laplacian }
\end{center}

 \bigskip\smallskip

\centerline{   Huyuan Chen\footnote{chenhuyuan@yeah.net},\ \ \quad Long Chen\footnote{ chenlong77@126.com}}
\medskip\bigskip
{\footnotesize
  \centerline{  Department of Mathematics, Jiangxi Normal University, Nanchang,}
  \centerline{  Jiangxi 330022, PR China}

}
\medskip

\begin{abstract}
  Our aim in this article is to obtain the  limit of  counting function for the Dirichlet eigenvalues  involving  the $m$-order logarithmic Laplacian $\cL_{m}$ in a bounded Lipschitz domain  $\Omega\in\mathbb{R}^{N}, N\ge1$ and to derive also the lower bound, where $m(\geq 2)\!\in\N$ and the Fourier transform $\cF(\mathcal{L}_{m})(\xi)=(2\log |\xi|)^m$.
\end{abstract}
 \bigskip

\noindent
  \noindent {\small {\bf Key Words}:   The $m$-order Logarithmic Laplacian; Counting function;  Dirichlet eigenvalues.  }\vspace{1mm}

\noindent {\small {\bf MSC2010}:  	35J10,  35R11, 35A01. }
\vspace{1mm}
\hspace{.05in}
\medskip

\setcounter{equation}{0}
\section{Introduction}
Let integer $N\geq1$,  $\Omega$ be a  bounded domain in $\R^N$ with the Lipschitz boundary for $N\geq2$.
Our purpose of this paper is to study the spectral of $\cL_{m}$ in a bounded domain, i.e.
 \begin{equation}\label{eq 1.1}
 \left\{
  \begin{array}{lll}
 \cL_{m}u =\lambda u\quad  &{\rm in}\ \, \Omega,
 \\[2mm] \phantom{   \cL_{m}  }
    u=0\quad    &{\rm in}  \ \, \R^N\setminus \Omega,
  \end{array}
 \right.
 \end{equation}
 where
  $\cL_{m}$ with $m\in \N$ is the $m$-order logarithmic Laplacian with
its Fourier transformation $\cF(\cL_{m})(\xi)=(2\log \xi)^m.$

In recent years, as a prototype model integro-differential operators, the fractional Laplacian have been studied in various aspects,   see e.g.
\cite{CS0,CS1,KM,MN,RS,S07} and the references therein.
Here the fractional Laplacian  $(-\Delta)^\s$ is defined as a singular integral operator
 \begin{equation}\label{fl 1}
 (-\Delta)^\s  u(x)=c_{N,\s} \lim_{\epsilon\to0^+} \int_{\R^N\setminus B_\epsilon(x) }\frac{u(x)-
u(y)}{|x-y|^{N+2\s}}  dy
\end{equation}
with a normalized constant  $c_{N,\s}=2^{2\s}\pi^{-\frac N2}\s\frac{\Gamma(\frac{N+2\s}2)}{\Gamma(1-\s)}.$
   An interesting limiting
property proposed in \cite{EGE} is the following
 \begin{equation} \label{eq:limit-behaviour}
\lim_{\s\to1^-}(-\Delta)^\s u(x)=-\Delta u(x)
\quad\  {\rm and}\quad\ \lim_{\s\to0^+}  (-\Delta)^\s  u(x) = u(x)\quad\  {\rm for} \ \, u\in C^2_c(\R^N).
 \end{equation}
    The author \cite{CW0} extended the expansion (\ref{eq:limit-behaviour}) to
    for $u \in C^2_c(\R^N)$ and  $x \in \R^N$,
 \begin{equation} \label{eq:limit-behaviour1}
(-\Delta)^\s  u(x) = u(x) + s \cL_1 u (x) + o(\s) \quad{\rm  as}\ \,  \s\to0^+,
 \end{equation}
where, formally, the operator $\cL_1$ is   the 1-orde {\em logarithmic Laplacian}. Moreover,  qualitative properties of $\cL_1$ are studied there.  Here we mention \cite{CDK,FJW,HS} on some related limiting properties of fractional problems.

Very recently, the first author of this paper in \cite{C0} gave a complete expansion by involving high order   logarithmic Laplacian $\cL_{m}$.
Now we introduce the  $m$-order logarithmic Laplacian $\cL_{m}$.  Before this,   we need to provide
some   notations:
 \begin{equation}\label{1.1-k0}
\kappa_1(s)  = 2^{-2s}\pi^{-\frac N2}  \frac{\Gamma(\frac{N-2s}{2})}{\Gamma(1+s)}
\quad {\rm and}\quad
\kappa_2(s)  = 2^{-2s}  \frac{\Gamma(\frac{N-2s}{2})}{\Gamma(\frac N2)\, \Gamma(1+s)}.
\end{equation}
Note that $\kappa_i\in C^\infty$ in $(-\min\{\frac{N}{2},1\}, \min\{\frac{N}{2},1\})$   and we use the notation $\kappa_i^{(m)}(0)$  the $m$-th order  derivative for $i=1,2$.

Given an integer $m\geq1$, the $m$ order logarithmic Laplacian  $\cL_m$ is an  integro-differential operator defined by
\begin{equation} \label{m 1.1}
\cL_m \zeta(x)=    \sum^{m}_{j=0}  \alpha_j \cK_{j} \zeta(x) \quad {\rm for}\  \, \zeta\in C^{0,\beta}_c(\R^N),
\end{equation}
where $\beta\in(0,1)$,
 $$\alpha_0=(-1)^{m}  \kappa_2^{(m)}(0), \qquad  \alpha_j=m (-1)^{m+j} \binom{m-1}{j-1} \kappa_1^{(m-j)}(0),\ \ j=1,2,\cdots,m  $$
and
\begin{align}\label{representation-main}
\cK_0  \zeta=\zeta, \qquad   \cK_j  \zeta (x)   =\int_{\R^N} \frac{\zeta(x)1_{B_1(x)}(y)-\zeta(y)}{|x-y|^N} (-2\log |x-y|)^{j-1} \,dy \ \ {\rm for\ integer }\ j\geq 1.
\end{align}
   With help of $\cL_m$, the expansion of fractional Laplacian with respect to the order could be written that   for $u \in C^2_c(\R^N)$ and  $x \in \R^N$,
$$
(-\Delta)^\s  u(x) = u(x) +\sum^n_{m=1}\frac{s^m}{m!} \cL_m u (x) + o(\s^n) \quad{\rm  as}\ \,  \s\to0^+.
$$
From this expression,  $\cL_m$ could be also viewed as the $m$-order derivative of fractional Laplacian, i.e.  $\mathcal{L}_{m}=\frac{d^{m}}{ds^{m}}(-\Delta)^{s}\big|_{s=0^{+}}$.

\smallskip

 As one of most important of properties, eigenvalues for differential operators has been studied extensively and deeply. {\it Denote by $\cN(\cL, \Omega,\mu)$ the number of Dirichlet eigenvalues of operator $\cL$ in $\Omega$, counted according their multiplicity, not exceeding the value $\mu$.}  In 1912, Weyl in \cite{W} shows that  the $k$-th eigenvalue $\mu_k(\Omega)$ of Dirichlet problem
 \begin{equation}\label{eq 1.1-lap}
 -\Delta u=\mu  u\quad     {\rm in}\ \   \Omega, \qquad
  u=0\quad   {\rm{on}}\  \   \partial \Omega
  \end{equation}
 has the asymptotic behavior
$$\cN(-\Delta, \Omega,\mu) \sim   b_N |\Omega| \mu^{\frac N2}\quad{\rm as}\ \,  \mu\to+\infty,$$ 
and
$$\mu_k(\Omega)\sim  c_N(k|\Omega|)^{\frac2N}\quad{\rm as}\ \,   k\to+\infty, $$
where $b_N=(2\pi)^{-\frac{N}{2}} |B_1|$ and $ c_N=(2\pi)^2|B_1|^{-\frac{2}{N}}$. 

 In 1960  P\'olya in \cite{P} gave an  equivalent   bound
\begin{equation}\label{p-conj}
\mu_k(\Omega)\geq c_N\left(\frac{k}{|\Omega|}\right)^{\frac{2}N},\quad \cN(-\Delta, \Omega,\mu) \leq  b_N |\Omega| \mu^{\frac N2}  
\end{equation}
for  any "plane-covering domain" $D$ in $\R^2$, this proof  works in dimension $N\geq 3$.
    From then,    it attracts great attentions for the  asymtoptics of the eigenvalues in various aspect. The authors in \cite{Ro,L} independently  proved first lower bound in (\ref{p-conj}) with a positive constant $c$ for general bounded domain;    Li-Yau \cite{LY} improved the value of
 the constant $C$ obtaining $C=\frac{N}{N+2} C_N$.   We also refer to   \cite{CgW,K,CgY} for generalized topics,   \cite{CL,HY}  for degenerate elliptic operators and  \cite{BCH,CS,CV0,YY}  for the fractional Laplacian, where 
  the lower bound of the eigenvalues reads
 \begin{equation}\label{p-conjf}
\mu_{s,k}(\Omega)\geq \frac{N}{N+2s}C_N\left(\frac{k}{|\Omega|}\right)^{\frac{2s}N}.
\end{equation}
Moreover,  in the fractional case,  \cite{Ge} extends the
limit of the counting function of the fractional Laplacian
$$\cN((-\Delta)^s, \Omega,\mu) \sim   b_N |\Omega| \mu^{\frac N{2s}}\quad {\rm as}\ \, \mu\to+\infty. $$

Involving the logarithmic Laplacian with the order $m=1$, the basic properties of Dirichlet eigenvalues in bounded domain  are investigated \cite{CV0,FJW},   the limit of counting function for the Dirichlet eigenvalues is derived in \cite{LW} 
$$
   \cN (\cL_1,\omega, \lambda)
\sim \frac{\omega_{_N}}{N(2\pi)^{N}}     |\Omega| e^{ \frac N2  \lambda}\quad {\rm as}\ \, \lambda\to+\infty. 
$$

Our aim of this article is to give  the limit of counting function and the lower bound of eigenvalues for high order logarithmic Laplacian with $m\geq2$.  For this,  we first denote by $\mathbb{H}_{m,0}(\Omega)$   the closure of $C^{\infty}_{c}(\Omega)$ with respect to the norm
 $$  \|u\| _{m}:=\sqrt{\int_{\mathbb{R}^{N}}(\log (e+|\xi|))^{m}|\hat{u}(\xi)|^{2}d\xi}.$$
It is shown in \cite{C0} that the embedding $\cH_{m,0}(\Omega)\hookrightarrow L^{2}(\Omega)$ is compact.
Note that   the corresponding  quadratic form of $\cL_m$ could be defined as
 $$\varphi \mapsto (\varphi,\varphi)_m :=\frac{1}{(2\pi)^{N}}\int_{\mathbb{R}^{N}}(2\log |\xi|)^{m}|\hat{\varphi}(\xi)|^{2}d\xi.$$
 It is shown that  $(\cdot,\cdot)_{m}$ defines a closed, symmetric and semibounded quadratic form with domain $\mathbb{H}_{m,0}(\Omega)$.
 Note that
 $$ \mathcal{D}(\cL_m)\subset \cH_{m,0}(\Omega) .$$
 The spectrum point  $(\lambda_{m,i},\phi_i)_{i\in \N}\in \R\times \cH_{m,0}(\Omega)$ of $\cL_{m}$ could  write as
 \begin{equation}\label{def id 1.1}
 (\phi_{m,k},v)_m=\lambda_{m,i} \int_{\Omega} \phi_{m,k}  v dx\quad {\rm for\ any}\  v\in \cH_{m,0}(\Omega).
 \end{equation}
 From \cite[Theorem 1.3]{C0} problem  (\ref{eq 1.1}) has a sequence of Dirichlet eigenvalues
 $$
 -\infty<\lambda_{m, 1}(\Omega)\leq \lambda_{m,2}(\Omega)\leq\cdots\leq\lambda_{m,k}(\Omega)\leq \cdots \quad \to+\infty\ \
 {\rm as}\ \, k\to+\infty  $$
 and  the corresponding eigenfunction $\phi_{m,k}\in \cH_{m,0}(\Omega)$ such that $\|\phi_{m,k}\|_{L^2(\Omega)}=1$.
 \begin{remark}
 From the Fourier transform $\cF(\cL_{m})(\xi)=(2\log \xi)^m$, we have the relationship that
 $$\cL_m=\underbrace{\cL_1\circ\cdots\circ \cL_1}_{m}.  $$
Letting for $(\lambda_{1,k},\phi_{1,k})_{k\in\N}$ be the Dirichlet eigenvalues and related eigenfunctions of $\cL_1$, then it fails to get
\begin{align*}
\cL_2\phi_{1,k} =  \cL_1 \circ (\cL_1 \phi_{1,k})  = \lambda_{1,k}^2 \phi_{1,k} \quad {\rm in}\ \, \Omega
\end{align*}
because although $\cL_1 \phi_{1,k}=\lambda_{1,k}\phi_{1,k}$ in $\Omega$ but $\cL_1 \phi_{1,k}\not=0$ in $\Omega^c$, because of the non-local property of $\cL_1$.   Therefore,  the set of the Dirichlet eigenvalues  of $\cL_m$  is different from  $\{(\lambda_{1,k}^m(\Omega)\}_{k\in\N}$.

 \end{remark}

 \smallskip

 For simplicity,   we use the notation of $\cN_{m}(\lambda)=\cN(\cL_m,\Omega, \lambda)$, i.e.
  \begin{align}\label{cont m 2.1}
 \cN_{m}(\lambda)=\# \big\{ k:\lambda_{m,k}(\Omega)< \lambda \big\},
 \end{align}
 which is  the counting function for  the Dirichlet eigenvalues  of the operator $\mathcal{L}_{m}$. 
Our main result  on the  limit of counting function  for the Dirichlet eigenvalues $\{\lambda_{m,k}(\Omega)\}_{k\in\N}$.

 \begin{theorem}\label{teo 1}
 Let $m\geq2$,  $\lambda_{m,k}(\Omega)$ be the $k^{th}$ eigenvalue of the operator $\mathcal{L}_{m}$ in $\Omega$ under the zero Dirichlet boundary condition, then
   \begin{align}\label{limt 4.0}
 \lim_{\lambda\to+\infty} e^{-\frac N2 \lambda^{\frac{1}{m}}} \cN_{m}(\lambda)
 =\frac{  T_N |\Omega|}{N} 
  \end{align}
  and
  \begin{align}\label{limt 4.0-+}
  \lim_{\lambda\to\infty}    \lambda^{-\frac{m-1}{m}} 
  e^{-\frac N2 \lambda^{\frac{1}{m}}} \sum_{k}(\lambda-\lambda_{m,k}(\Omega))_{+}  = \frac{2m}{N^2} T_N  |\Omega|,
  \end{align}
where $T_N=\frac{\omega_{_N}}{(2\pi)^{N}}$ and  $a_+=\max\{0,a\}$.
  \end{theorem}


We provide asymptotic formulas for the counting function (\ref{limt 4.0}) and (\ref{limt 4.0-+}), which could lead to  the Weyl's limit if we have the estimates of multiplicity of eigenvalues.   Here we want to mention that the Weyl's limit is derived   from the upper and lower bounds in \cite{CV0} for $m=1$ by the Li-Yau's lower bound and Kr\"oger's upper bound.  An very interesting point is that the limit of the sum of first $k$ eigenvalues coincides the Weyl's limit and  they don't depends on the volume of the domain.
 \smallskip

 For lower bound, we have the following estimates: 
 \begin{theorem}\label{teo 2}
 Let the integer $m,k\geq 2$,  $\lambda_{m,k}(\Omega)$ be the $k^{\rm th}$ Dirichlet eigenvalue of the operator $\mathcal{L}_{m}$,
 $$a_m= \frac{N^{m+1}    }{ 2^m  m} \frac{(2\pi)^{N}}{ \omega_{_N}},\quad  b_m=  \frac{2^m}{N^m}  (m-1)^m +\frac1{N^m} $$
 and  $c_m=a_mb_m$.

$(i)$ When
  either  $m$ is odd and $\lambda_{m,1}(\Omega)\geq 0$ or  $m$ is even,
then  \begin{align*}
\lambda_{m,k}(\Omega)  \geq  &\big(\frac{2}{N}\big)^m   \min\Big\{  \big(  \frac{c_m}{e^{\tau_0}} \big)^{\frac m{m-1}}      ( \frac{ k }{|\Omega|})^{\frac m{m-1}}, \nonumber\\[2mm]
 &\qquad\qquad\quad\ \Big(\log \big(c_m  \frac{k }{|\Omega|}    + \tau_0 ^{m-1}  e^{\tau_0}\big)  - (m-1)\log \big(\log (c_m  \frac{k }{|\Omega|} + e)  \big) -\log2 \Big)^m    \Big\}   - b_m,
  \end{align*}
  where      $\tau_0 > 1$ is only dependent of $m$ and $\omega_{_N}$ is  the volume of the unit sphere in $\mathbb{R}^{N}$. \smallskip

$(ii)$ When   $m$ is odd and $\lambda_{m,1}(\Omega)<0$,  then
  \begin{align*}
\lambda_{m,k}(\Omega)  \geq  \big(\frac{2}{N}\big)^m   \min\Big\{ &\big(  \frac{a_m}{e^{\tau_0}} \big)^{\frac m{m-1}}   \big(\frac{ P_m (\Omega) }{  |\Omega|} \big)^{\frac m{m-1}}   k^{\frac m{m-1}}, \\
 &\Big(\log (a_m \frac{P_m(\Omega)  }{|\Omega|}  k + \tau_0 ^{m-1}  e^{\tau_0})  - (m-1)\log \big(2\log (a_m  \frac{P_m(\Omega) }{|\Omega|}  k + e)   \big)   \Big)^m    \Big\}   - P_m(\Omega),
  \end{align*}
where $P_m(\Omega)=- \lambda_{m,1}(\Omega)   +b_m$.

  \end{theorem}

  Our main difficulty comes from $m\geq 2$ in the derivation of the  limits of counting function and the lower bound.  For $m=1$, there is a special scaling property
$$\lambda_{1,k}(R\Omega)=\lambda_{1,k}(\Omega)-\log R,$$
because of the property $\log (R|\xi|)=\ln(|\xi|)+\log R$, where $R\Omega=\{Rx:\, x\in\Omega\}$.
However, the above property fails for $m\geq2$ and we provide a related rescaling estimates.


 Due to  the inhomogeneity of the expression of the  Fourier symbol of $(2\log(|\!\cdot\!|)^m$,    the essential tool for the lower bound is the estimates of  the count  functionfunction.  It is worth noting that our lower bound $\lambda_{m,k}$ is based on the value $\lambda_{m,1}$ and we address some basic properties of  $\lambda_{m,1}$ in Section 2.   In Section 3, we build the basic estimates for  the count function of the eigenvalues and the   limit of counting functions.
 Section 4 is devoted to    lower bound of $\lambda_{m,k}(\Omega)$.

\setcounter{equation}{0}
 \section{ Preliminary  properties}

  \subsection{ Recaling property }

We start the section from the recalling property.
\begin{proposition}\label{lm rec}
Let $m\geq 3$ be {\bf odd},  $\Omega$ be a Lipschitz  bounded open set and   $\lambda_{m,k}(\Omega)$ be the $k^{th}$ Dirichlet eigenvalue of the operator $\mathcal{L}_{m}$, for $k\geq 2$. For $R>1$, denote
$$\Omega_{_R}=\{Rx: x\in\Omega\},$$
then
$$2^{-m}  \lambda_{m,k}(\Omega) - 4^m (\log R)^m -\frac{4^m-1}{(2\pi)^{N} }C_m  |\Omega|  \leq \lambda_{m,k}(\Omega_{_R})\leq  \lambda_{m,k}(\Omega)-2 (\log R)^m, $$
where $C_m= \int_{B_1} \big|\log|\xi|\big|^m d\xi$.
\end{proposition}
{\bf Proof. } Due to $C_c^{\infty}(\Omega) \subset \cH(\Omega)$ is dense, we set $\phi_{m,k} \in C_c^{\infty}(\Omega)$ and $\phi_{R,m,k} (x)= R^{-\frac{N}{2}}  \phi_{m,k}(\frac{x}{R})\in C_c^{\infty}(\Omega)$, \\
then $\| \phi_{R,m,k}\|_{L^2(\R^N)}=\| \phi_{m,k}\|_{L^2(\R^N)}$, and
\begin{align*}
\hat{\phi}_{R,m,k} (\zeta)&=\int_{\R^N} e^{-i\zeta \cdot x} R^{-\frac{N}{2}}  \phi_{m,k}(\frac{x}{R}) dx
\\[1mm]&=\int_{\R^N} e^{-iR \zeta \cdot y} R^{-\frac{N}{2}}  \phi_{m,k}(y) d(Ry)
\\[1mm]&=R^{\frac{N}{2}} \int_{\R^N} e^{-i R \zeta \cdot y}   \phi_{m,k}(y) dy
\\[1mm]&=R^{\frac{N}{2}} \hat{\phi}_{m,k} (R\zeta),
 \end{align*}
where $y=\frac{x}{R}$, then we obtain that
\begin{align*}
\lambda_{m,k}(\Omega_{_R}) &=\lambda_{m,k}(\Omega_{_R}) \int_{\R^N} |\phi_{R,m,k} (x)|^2 dx
\\[1mm]&=\int_{\R^N} \cL_m (\phi_{R,m,k} (x)) \phi_{R,m,k} (x) dx
\\[1mm]&=\frac1{(2\pi)^N} \int_{\R^N} (2 \log|\zeta|)^m  |\hat{\phi}_{R,m,k} (\zeta)|^2 d\zeta
\\[1mm]&=\frac{2^m R^N}{(2\pi)^N} \int_{\R^N} ( \log|\zeta|)^m  |\hat{\phi}_{m,k} (R \zeta)|^2 d\zeta .
  \end{align*}

On the one hand, for $a,b>0$, $m\geq 3$ odd, one has
 $$(a+b)^m\geq  a^m+b^m \quad {\rm and}\quad (a-b)^m\leq  a^m-2^{1-m}b^m$$
thus, letting $\xi= R \zeta$, direct computation shows that
\begin{align*}
\lambda_{m,k}(\Omega_{_R})
&=\frac{2^m}{(2\pi)^{N} } \int_{\R^N} \big(\log|\xi| -\log R\big)^m |\hat{\phi}_{m,k}(\xi)|^2 d\xi
\\[1mm]&=\frac{2^m}{(2\pi)^{N} } \int_{B_1} \big(\log|\xi| -\log R\big)^m |\hat{\phi}_{m,k}|^2 d\xi+\frac{2^m}{(2\pi)^{N} } \int_{\R^N\setminus B_1} \big(\log|\xi| -\log R\big)^m  |\hat{\phi}_{m,k}|^2 d\xi
 \\[1mm]&\leq  -\frac{2^m}{(2\pi)^{N} } \int_{B_1} \big(\big|\log|\xi|\big|^m +(\log R)^m\big)  |\hat{\phi}_{m,k}|^2 d\xi
 \\[1mm]&\quad +\frac{2^m}{(2\pi)^{N} } \int_{\R^N\setminus B_1} \big(\big|\log|\xi|\big|^m -2^{1-m}(\log R)^m\big)  |\hat{\phi}_{m,k}|^2 d\xi
  \\[1mm]&\leq \frac{2^m}{(2\pi)^{N} } \int_{\R^N}  \big|\log|\xi|\big|^m |\hat{\phi}_{m,k}|^2 d\xi  -\frac{2^m}{(2\pi)^{N} } \int_{B_1}  \big|\log|\xi|\big|^m |\hat{\phi}_{m,k}|^2 d\xi
  \\[1mm]&\quad  -2^{1-m} (\log R)^m \frac{2^m}{(2\pi)^{N} } \int_{\R^N}    |\hat{\phi}_{m,k}|^2 d\xi
  \\[1mm]&\leq  \lambda_{m,k}(\Omega) -2   (\log R)^m.
 \end{align*}

 On the other hand,
 for $a,b>0$, $m\geq 3$ odd, there holds
 $$(a+b)^m\leq 2^m( a^m+b^m) \quad {\rm and}\quad (a-b)^m\geq 2^{-m} (a^m-2^m b^m), $$
 then
 \begin{align*}
\lambda_{m,k}(\Omega_{_R})
 &=\frac{2^m}{(2\pi)^{N} } \int_{B_1} \big(\log|\xi| -\log R\big)^m |\hat{\phi}_{m,k}|^2 d\xi+\frac{2^m}{(2\pi)^{N} } \int_{\R^N\setminus B_1} \big(\log|\xi| -\log R\big)^m  |\hat{\phi}_{m,k}|^2 d\xi
 \\[1mm]&\geq  -\frac{4^m}{(2\pi)^{N} } \int_{B_1} \big(\big|\log|\xi|\big|^m +(\log R)^m\big)  |\hat{\phi}_{m,k}|^2 d\xi+\frac{1}{(2\pi)^{N} } \int_{\R^N\setminus B_1} \big(\big|\log|\xi|\big|^m -2^m (\log R)^m\big)  |\hat{\phi}_{m,k}|^2 d\xi
 \\[1mm]&\geq   \frac{1}{(2\pi)^{N} } \int_{\R^N}  \big(\log|\xi|\big)^m   |\hat{\phi}_{m,k}|^2 d\xi+\frac{4^m-1}{(2\pi)^{N} } \int_{B_1}  \big(\log|\xi|\big)^m   |\hat{\phi}_{m,k}|^2 d\xi-\frac{4^m}{(2\pi)^{N} }(\log R)^m \int_{\R^N}      |\hat{\phi}_{m,k}|^2 d\xi
  \\[1mm]&\geq 2^{-m}\lambda_{m,k}(\Omega) -2^m (2\log R)^m-\frac{4^m-1}{(2\pi)^{N} }C_m  |\Omega|,
 \end{align*}
 where
  \begin{align*}
\Big|  \int_{B_1}  \big(\log|\xi|\big)^m   |\hat{\phi}_{m,k}|^2 d\xi\Big|&\leq C_m |\hat{\phi}_{m,k}|_{L^\infty}^2     \\& \leq  C_m    | \phi_{m,k}|_{L^1}^2\leq C_m  |\Omega| | \phi_{m,k}|_{L^2}^2=C_m  |\Omega|.
  \end{align*}
  We complete the proof. \hfill$\Box$

\begin{corollary}\label{cr 2.1}
$(i)$ Let   $m\geq 3$ odd.  If $\Omega$ is small enough,  then
$\lambda_{m,1}(\Omega)>0$.  \smallskip

$(ii)$ Let   $m\geq 3$ odd,  $O$ be  a Lipschitz domain with fixed volume and $\Omega= rO$,  then
$\lambda_{m,1}(\Omega)<0$ for $r>1$ large enough.

\end{corollary}
{\bf Proof. }  It infers by \cite[Theorem 1.3$(v)$]{C0} that  $\lambda_{m,1}(\Omega) \geq0$   when $\Omega\subset B_{r_0}$ for some $r_0>0$ small.  We apply  Proposition \ref{lm rec} to obtain that
 $$\lambda_{m,1}(\Omega )=\lambda_{m,1}(rO)\leq\lambda_{m,1}(O) -2 (\log r)^m,$$
 then, for $r>1$ large enough, we have
$\lambda_{m,1}(\Omega )<0. $
 \hfill$\Box$

\subsection{Lower bound for $\lambda_{m,1}(B_{r_0})$ }
Let $r_0\in(0,1]$ be  such that 
 the   sum of the kernel 
 \begin{equation}\label{hhh}
  \sum^m_{j=1} \alpha_j  (-2\ln t)^{j-1} >0\quad {\rm for}\ \, 0<t<{r_0}. 
  \end{equation}

\begin{proposition}
Let integer $N\geq 2$,  odd integer $m\geq 3$ and $r_0\in(0,1)$,
 then we have
\begin{align}
\lambda_{m,1}(B_{r_0})&\geq 2^m (\log (2\sqrt{N+2})-\log r_0)^m \nonumber
\\[2mm]&\quad -\frac{2^m   (2\sqrt{N+2})^N  \omega_{_N}^2}{(2\pi)^{2N} } \sum\limits_{j=1}^{m}
\frac{(-1)^{j+1} }{N^{j+2}}  \frac{m!}{(m-j)!} (\log (2\sqrt{N+2})-\log r_0)^{m-j}.\label{abb-1}
\end{align}
 \end{proposition}
{\bf Proof. }Let $v\in L^2(B_{r_0})$ be a radial function such that  with $\| v\|_{L^2(B_{r_0})}=1$ and recall the Bessel function $J_l(t)=(\frac t2)^l  \sum\limits_{j=0}^{\infty} \frac{(-1)^j}{j! \Gamma(j+l+1)} (\frac t2)^{2j}$.
For $l\geq \sqrt{3}-2$ and $0\leq t\leq 2\sqrt{2(l+2)}$,  it follows by the appendix of \cite{LW}
$$|J_l(t)|\leq \frac{t^l}{2^l \Gamma(l+1)}.$$
Taking $l=\frac N2  -1$,  we have that
$$|J_{\frac N2  -1} (t)|\leq \frac{t^{\frac N2  -1}}{2^{\frac N2  -1} \Gamma(\frac N2)}=\frac{\omega_{_N}}{(2\pi)^{\frac N2}} t^{\frac N2  -1}.$$

Note that  $\hat{v}$ is also radial and
\begin{align}
|\hat{v}(\xi)|&=s^{1-\frac N2} \Big| \int_0^{r_0} v(r) J_{\frac N2  -1} (rs) r^{\frac N2} dr\Big|
\nonumber\\[1mm]&\leq s^{1-\frac N2} \Big(\int_0^{r_0}  r^{N-1} v^2(r) dr  \Big)^{\frac 12}  \Big(\int_0^{r_0} r J^2_{\frac N2  -1} (rs)  dr  \Big)^{\frac12}\label{xxx-1-0}
\\[1mm]&=\frac{s^{1-\frac N2}}{\omega_{_N}^{\frac12}}  \Big(s^{-2} \int_0^{r_0s} t J^2_{\frac N2  -1} (t)  dt  \Big)^{\frac12},\nonumber
 \end{align}
 where  $s=|\xi|$ and $t=sr$.
 Then
 $$\omega_{_N} |\hat{v}(\xi)|^2 \leq s^{-N} \int_0^{r_0s}  t J^2_{\frac N2  -1} (t)  dt.$$

Now we take $v$  the  first eigenfunction of (\ref{eq 1.1}) corresponding to $\lambda_{m,1}$ in $\Omega=B_{r_0}$,
which is radial symmetric and decreasing with respect to $|x|$ by \cite[Theorem 1.3 $(v)$]{C0}.  Then  for every $\lambda \in \R$,
 \begin{align*}
 (2\pi)^N \Big(\lambda-\lambda_{m,1}(B_{r_0})  \Big)
        &=\int_{\R^N} (\lambda-(2\log|\xi|)^m) |\hat{v}(\xi)|^2 d\xi
 \\[1mm]&\leq \int_{\R^N} (\lambda-(2\log|\xi|)^m)_+ |\hat{v}(\xi)|^2 d\xi
 \\[1mm]&= \omega_{_N} \int_0^{e^{\frac12 \lambda^{\frac1m}}} (\lambda-(2\log s)^m)   |\hat{v}(s)|^2   s^{N-1} ds
 \\[1mm]&\leq  \int_0^{e^{\frac12 \lambda^{\frac1m}}} \frac{(\lambda-(2\log s)^m)}{s} \int_0^{r_0 s}  t J^2_{\frac N2  -1} (t)  dtds
 \\[1mm]&=\int_0^{r_0 e^{\frac12 \lambda^{\frac1m}}} t J^2_{\frac N2  -1} (t)  \int_{\frac{t}{r_0}}^{e^{\frac12 \lambda^{\frac1m}}} \frac{(\lambda-(2\log s)^m)}{s} dsdt
 \\[1mm]&=\int_0^{r_0 e^{\frac12 \lambda^{\frac1m}}} t J^2_{\frac N2  -1} (t)  \big(\frac m{2m+2} \lambda^{\frac{m+1}{m}} - \lambda \log \frac{t}{r_0} +\frac{2^m}{m+1} (\log \frac{t}{r_0})^{m+1} \big)dt
 \\[1mm]&=\int_0^1 r_0^2 e^{ \lambda^{\frac1m}} t J^2_{\frac N2  -1} (r_0 e^{\frac12 \lambda^{\frac1m}}t)  \big(-\frac {1}{2m+2} \lambda^{\frac{m+1}{m}} - \lambda \log t +\frac{2^m}{m+1} (\log t+\frac12 \lambda^{\frac1m})^{m+1} \big)dt
 \\[1mm]&\leq \frac{\omega_{_N}^2 r_0^N   e^{\frac N2 \lambda^{\frac1m}}}{(2\pi)^N} \int_0^1 \big(-\frac {1}{2m+2} \lambda^{\frac{m+1}{m}} - \lambda \log t +\frac{2^m}{m+1} (\log t+\frac12 \lambda^{\frac1m})^{m+1} \big) t^{N-1} dt
 \\[1mm]&=\frac{\omega_{_N}^2 r_0^N   e^{\frac N2 \lambda^{\frac1m}}}{(2\pi)^N N(m+1)} \int_0^1 \sum\limits_{j=2}^{m+1} \frac{2^{j-1}  (m+1)!}{(m+1-j)! \, j!} \lambda^{\frac{m+1-j}{m}} (\log t)^j  dt^N
 \\[1mm]&=\frac{\omega_{_N}^2  r_0^N  e^{\frac N2 \lambda^{\frac1m}}}{(2\pi)^N } \sum\limits_{j=2}^{m+1} (-1)^j   \frac{2^{j-1}}{N^{j+1}}  \frac{m!}{(m+1-j)!} \lambda^{\frac{m+1-j}{m}}.
  \end{align*}
 Hence
 \begin{align}\label{abb-0}
 \lambda_{m,1}(B_{r_0})\geq \lambda - \frac{\omega_{_N}^2 r_0^N   e^{\frac N2 \lambda^{\frac1m}}}{(2\pi)^{2N} } \sum\limits_{j=1}^{m} (-1)^{j+1}   \frac{2^j}{N^{j+2}}  \frac{m!}{(m-j)!} \lambda^{\frac{m-j}{m}}.  
 \end{align}
If $N\geq2$, $r_0  e^{\frac12 \lambda^{\frac1m}}\leq 2\sqrt{N+2}$, then
$$\lambda\geq 2^m (\log(2\sqrt{N+2})-\log r_0)^m.$$
Taking $\lambda= 2^m (\log (2\sqrt{N+2})-\log r_0)^m$,  (\ref{abb-0}) implies (\ref{abb-1}). 
  \hfill$\Box$

\setcounter{equation}{0}
\section{ Limit for counting function}

\subsection{bounds for counting function}


The upper bound of $\cN_{m}$ states as following:
    \begin{lemma}\label{pr 2.1}
Let  $\cN_{m}$ be given in (\ref{cont m 2.1}), then for any $\eta > \lambda >  \lambda_{m,k}(\Omega)$
 $$\cN_{m}(\lambda) \leq
  T_N|\Omega| \frac{1}{\eta-\lambda}  e^{\frac N2 \eta^{\frac{1}{m}}}
 \Big( \sum \limits^{m}_{j=1} (-1)^{j+1} \frac{A_{m,j}}{N^{j+1}} \eta^{\frac{m-j}{m}} \Big),$$
  where
 \begin{align}\label{aaa 2.1}
 A_{m,j}=2^j\frac{m!}{(m-j)!}.
 \end{align}

  \end{lemma}
   {\bf Proof. }  We first prove that
 \begin{align}\label{upp 2.1}
 \sum_{k}(\lambda-\lambda_{m,k}(\Omega))_{+} \leq \frac{|\Omega|\omega_{_N}}{(2\pi)^{N}}e^{\frac N2 \lambda^{\frac{1}{m}}}
 \Big( \sum^{m}_{j=1} \frac{ (-1)^{j+1} A_{m,j}}{N^{j+1}} \lambda^{\frac{m-j}{m}} \Big),
 \end{align}
 where we recall $a_+=\max\{0,a\}$. \smallskip

 Note that the eigenfunction $\phi_{m,k}$ has the zero data outside of $\Omega$, then  the Fourier transform implies that
 \begin{align*}
 \sum_{k}(\lambda-\lambda_{m,k}(\Omega))_{+}&=\sum_{k}\big(\lambda(\phi_{m,k},\phi_{m,k})_{L^{2}(\Omega)}
 -(\mathcal{L}_{m}\phi_{m,k},\phi_{m,k})_{L^{2}(\Omega)}\big)_{+}\\
 &=\frac{1}{(2\pi)^{N}}\sum_{k}\Big(\int_{\mathbb{R}^{N}}\big(\lambda-(2\log|\xi|)^{m}\big)|\hat{\phi}_{m,k}(\xi)|^{2}d\xi\Big)_{+}\\
 &\leq \frac{1}{(2\pi)^{N}}\sum \limits_{k}\int_{\mathbb{R}^{N}}\big(\lambda-(2\log|\xi|)^{m}\big)_{+}|\hat{\phi}_{m,k}(\xi)|^{2}d\xi\\
 &=\frac{1}{(2\pi)^{N}}\int_{\mathbb{R}^{N}}\big(\lambda-(2\log|\xi|)^{m}\big)_{+} \sum\limits_{k} |\hat{\phi}_{m,k}(\xi)|^{2}d\xi.
 \end{align*}
 Since $\{ \phi_{m,k} \}$ is an orthonormal basis in $L^{2}(\Omega)$ and let $e_{\xi}(x)=e^{-{\rm i} x\cdot\xi }$, then by Pythagorean theorem:
 $$\sum\limits_{k}|\hat{\phi}_{m,k}(\xi)|^{2}=\sum\limits_{k}|(e_{\xi},\phi_{m,k})_{L^{2}(\Omega)}|^{2}
 =\|e_{\xi}\|^{2}_{L^{2}(\Omega)}=|\Omega|, $$
 letting $r_m=e^{\frac12 \lambda^{\frac1m}}$, and we have that
 \begin{align*}
 \sum\limits_{k}(\lambda-\lambda_{m,k}(\Omega))_{+}
 &\leq \frac{1}{(2\pi)^{N}}|\Omega|\int_{\mathbb{R}^{N}}(\lambda-(2\log|\xi|)^{m})_{+}d\xi\\
 &=\frac{1}{(2\pi)^{N}}|\Omega|\int_{B_{r_m}(0)}\Big(\lambda-(2\log|\xi|)^{m}\Big)d\xi\\
 &= T_N |\Omega| e^{\frac N2 \lambda^{\frac{1}{m}}}
 \Big(\sum \limits^{m}_{j=1} (-1)^{j+1} \frac{A_{m,j}}{N^{j+1}} \lambda^{\frac{m-j}{m}}\Big),
 \end{align*}
where we used the fact that
 \begin{align*}
 \int_{B_{r_m}}(2\log |\xi|)^m d\xi&=2^m \omega_{_N} \int_0^{r_m} (\log r)^m r^{N-1} dr
 \\& = 2^m \omega_{_N} r_m^N \Big(\frac1N (\log r_m)^m+  \sum^m_{j=1} \frac{(-1)^j}{N^{j+1}} \frac{m!}{(m-j)!}  (\log r_m)^{m-j}\Big).
  \end{align*}

 Thus, we obtain that
 $$\cN_{m}(\lambda) \leq \frac{1}{\eta-\lambda} \sum\limits_{k}(\eta-\lambda_{m,k}(\Omega))_{+} \leq
 T_N |\Omega|   \frac{1}{\eta-\lambda}  e^{\frac N2\eta^{\frac{1}{m}}}
 \Big( \sum \limits^{m}_{j=1} (-1)^{j+1} \frac{A_{m,j}}{N^{j+1}} \eta^{\frac{m-j}{m}} \Big).$$
 We complete the proof.
  \begin{remark}
 Under  the estimate (\ref{upp 2.1}), we can see that  
  $$ \lim_{\lambda \to +\infty} \frac{\sum\limits^{m}_{j=1} (-1)^{j+1} \frac{A_{m,j}}{N^{j+1}} \lambda^{\frac{m-j}{m}}}{\frac{A_{m,1}}{N^2} \lambda^{\frac{m-1}{m}}}=1,$$
  then
  $$\limsup_{\lambda \to +\infty} \frac{\sum\limits_{k}(\lambda-\lambda_{m,k}(\Omega))_{+}}{e^{\frac N2 \lambda^{\frac{1}{m}}} \frac{A_{m,1}}{N^2} \lambda^{\frac{m-1}{m}}}
  \leq
  \limsup_{\lambda \to +\infty} \frac{\frac{|\Omega|\omega_{_N}}{(2\pi)^{N}} e^{\frac N2 \lambda^{\frac{1}{m}}}  \Big( \sum^{m}_{j=1} \frac{ (-1)^{j+1} A_{m,j}}{N^{j+1}} \lambda^{\frac{m-j}{m}} \Big)}{e^{\frac N2 \lambda^{\frac{1}{m}}} \frac{A_{m,1}}{N^2} \lambda^{\frac{m-1}{m}}}
  \leq  T_N |\Omega|.$$
  then we have that
  \begin{align}\label{re2.5}
  \lim_{\lambda \to +\infty} \Big(  \frac{2m}{N^2} \lambda^{\frac{m-1}{m}}\Big)^{-1}
   e^{-\frac N2 \lambda^{\frac{1}{m}}}\sum_{k}(\lambda-\lambda_{m,k}(\Omega))_{+}  \leq   T_N |\Omega| .
  \end{align}
  \end{remark}

  \hfill$\Box$\medskip

Now we have the  asymptotic lower  bound for $(\lambda-\lambda_{m,k})_{+} $.

 \begin{lemma}\label{pr 2.2}
 For the eigenvalues of the problem (1.1) and any $\lambda\in \mathbb{R}$ we have
 \begin{align}\label{es 2.2}
 \liminf_{\lambda\to+\infty} \Big(  \frac{2m}{N^2} \lambda^{\frac{m-1}{m}}\Big)^{-1}
 e^{-\frac N2 \lambda^{\frac{1}{m}}}\sum_{k}(\lambda-\lambda_{m,k}(\Omega))_{+}  \geq   T_N |\Omega| .
 \end{align}

  \end{lemma}

  To order to prove Lemma \ref{pr 2.2},  we need  introduce the "coherent state"
 \begin{align}\label{coh 2.1}
 e_{\xi,y}(x)=e^{-i\xi \cdot x}\eta_\delta (x-y)\quad{\rm for}\ \,  \xi\in\R^N,\ x \in \Omega_{\delta},
 \end{align}
where  $\delta > 0$,
 $$\Omega_{\delta}=\big\{ x\in \Omega: \, {\rm dist}(x,\mathbb{R}^{N}\setminus \Omega)>\delta \big\} $$
 and
 $\eta_\delta:\R^N\to [0,+\infty) $  is a radially symmetric,  decreasing with respect to $|x|$,  function  in $C^{\infty}_0(\R^N)$  such that
   $\|\eta_\delta\|_{L^{2}(\mathbb{R}^{N})}=1$ with support in $B_{\frac\delta2}(0)$.
Note that  for any $M>0$,   there exists $c_M=\frac{(2\pi)^N \omega_{_N}}{N} (\frac{\delta}2)^N (1+|\frac{\delta}2|)^{M}>0$ such that
\begin{align*}
|\hat{\eta}_\delta(\xi)|^{2}
=\Big(\int_{\mathbb{R}^{N}} e^{-i \xi \cdot x} \,  {\eta}_\delta(x) dx \Big)^2
&=\Big(\int_{B_{\frac{\delta}2}} e^{-i \xi \cdot x}  \,  {\eta}_\delta(x) dx \Big)^2 \\
&\leq \int_{B_{\frac{\delta}2}} |e^{-i \xi \cdot x}|^2 dx  \, \int_{B_{\frac{\delta}2}} | {\eta}_\delta(x)|^2 dx \\
&=\frac{(2\pi)^N \omega_{_N}}{N} (\frac{\delta}2)^N \\
&\leq \frac{(2\pi)^N \omega_{_N}}{N} (\frac{\delta}2)^N (1+|\xi|)^{M}\\
&\leq \frac{(2\pi)^N \omega_{_N}}{N} (\frac{\delta}2)^N (1+|\frac{\delta}2|)^{M}
=c_M, \end{align*}

 thus for $\xi \in B_{\frac{\delta}2}(0)$, we have
 $$ |\hat{\eta}_\delta(\xi)|^{2} \leq \frac{(2\pi)^N \omega_{_N}}{N} (\frac{\delta}2)^N \leq c_M   (1+|\xi|)^{-M},$$

For given $M>0$, we let
 $$\beta_{m,j} =   T_N  c_{M}   \binom{m}{j}  \frac{1}{m-j+1}, \ \ j=0,1,2\cdots,m, $$
where   the binomial coefficient
    $$\binom{m}{i}=\frac{m!}{i!(m-i)!}$$
with  $i!$ being the $i$  factorial and $0!=1$. Particularly, $\beta_{m,m}=  T_N  c_{M}$.

    For $ \lambda> 2^m \beta_{m,m},$ we denote
\begin{align}\label{fun b1}
 f_\lambda (r)=\lambda-2^m (\log r)^m - 2^m\sum \limits^{m-1}_{j=1} \beta_{m,j}  (\log r)^{m-j}  r^{-j} - 2^m \beta_{m,m} r^{-m}  \quad{\rm for}\  r>0.
 \end{align}

 Note that $ f_\lambda (1)=\lambda - 2^m \beta_{m,m} > 0 $ for $ \lambda> 2^m \beta_{m,m}$,
 and
 $$\lim_{r\to 0^+}f_\lambda (r)=\lim_{r\to +\infty}f_\lambda (r)=-\infty,   $$
then   $f_\lambda $ has a minimum zero point $r_1=r_1(\lambda)\in(0,1)$ and a maximum zero point $r_2=r_2(\lambda)$ in $(1,+\infty)$.

Let
$$r_0=e^{\frac12 \lambda^{\frac1m}}>1, $$
which is the zero point of $g_{\lambda}(r):= \lambda (\log r)^{-m} - 2^m  \quad {\rm for}\ r>1$.
 Since $2^m\sum \limits^{m-1}_{j=1} \beta_{m,j}  (\log r)^{m-j}  r^{-j} >0$ for $r>1$,
 then
$$
g_{\lambda}(r) >  f_\lambda (r) (\log r)^{-m}\quad{\rm for}\ \, r>1,
$$
which implies
$
r_2< r_0 \ \, {\rm for}\ \lambda>0. $ \medskip

\noindent {\bf Proof of Lemma \ref{pr 2.2}. }  Note that $\|e_{\xi,y}\|_{L^{2}(\mathbb{R}^{N})}=1$. Using the properties of coherent states:
 $$\sum\limits_{k}(\lambda-\lambda_{m,k}(\Omega))_{+}\geq \frac{1}{(2\pi)^{N}}\int_{\mathbb{R}^{N}}\int_{\Omega_{\delta}}\big(e_{\xi,y},(\lambda-\mathcal{L}_{m})_{+} e_{\xi,y}\big)_{L^{2}(\Omega)} dy d\xi.$$
 where $e_{\xi,y}$ is given in (\ref{coh 2.1}).
 Since $t\mapsto(\lambda-t)_{+}$ is convex then applying Jensen's inequality to the spectral measure of $\mathcal{L}_{m}$ $$\sum\limits_{k}(\lambda-\lambda_{m,k}(\Omega))_{+}\geq \frac{1}{(2\pi)^{N}}\int_{\mathbb{R}^{N}}\int_{\Omega_{\delta}}\big(\lambda-(\mathcal{L}_{m}e_{\xi,y},e_{\xi,y})_{L^{2}(\Omega)}\big)_{+} dy d\xi.$$
 Next we consider the quadratic form
 \begin{align*}
 &\quad (\mathcal{L}_{m}e_{\xi,y},e_{\xi,y})_{L^{2}(\Omega)}
 \\[1mm]&= \frac{1}{(2\pi)^{N}}\int_{\mathbb{R}^{N}} (2\log|\gamma|)^{m} |\hat{e}_{\xi,y}(\gamma)|^{2} d\gamma
 \\[1mm]&= \frac{1}{(2\pi)^{N}}\int_{\mathbb{R}^{N}} (2\log|\gamma|)^{m} \overline{\hat{e}_{\xi,y} (\gamma)} \,\hat{e}_{\xi,y}(\gamma)dy d\gamma
 \\[1mm]&= \frac{2^m}{(2\pi)^{N}}\int_{\mathbb{R}^{N}}\Big(\int_{\mathbb{R}^{N}}e^{i\gamma\cdot x} \, e^{i\xi\cdot x} \eta_\delta(x-y)dx\Big)
 \Big(\int_{\mathbb{R}^{N}}e^{-i\gamma\cdot z }  e^{-i\xi\cdot z} \eta_\delta(z-y)dz\Big) ( \log|\gamma|)^{m} d\gamma
 \\[1mm] &= \frac{2^m}{(2\pi)^{N}} \int_{\mathbb{R}^{N}} \int_{\mathbb{R}^{N}} \int_{\mathbb{R}^{N}} e^{i(x-z)\cdot (\gamma+\xi)} \eta_\delta(x-y)  \eta_\delta(z-y)
 (\log|\gamma|)^{m} dx dz d\gamma
 \\[1mm] &= \frac{2^m}{(2\pi)^{N}} \int_{\mathbb{R}^{N}} \int_{\mathbb{R}^{N}} \int_{\mathbb{R}^{N}} e^{i(x-z)\cdot \rho} \eta_\delta(x-y) \eta_\delta(z-y) (  \log|\rho-\xi|)^{m} dx dz d\rho
 \\[1mm]
 &= \frac{2^m}{(2\pi)^{N}} \int_{\mathbb{R}^{N}} (\int_{\mathbb{R}^{N}} e^{i \rho \cdot (x-y)} \eta_\delta(x-y) dx) (\int_{\mathbb{R}^{N}} e^{-i\rho \cdot (z-y)} \eta_\delta(z-y) dz) \big[\log|\xi|+\log( \frac{|\rho-\xi|}{|\xi|})\big]^{m} dx dz d\rho
 \\[1mm]
 &= \frac{2^m}{(2\pi)^{N}} \int_{\mathbb{R}^{N}} |\hat{\eta_\delta}(\rho)|^{2} \Big[\sum \limits^{m}_{j=0}
 \binom{m}{j} (\log|\xi|)^{m-j} (\log( \frac{|\rho-\xi|}{|\xi|}))^{j} \Big] d\rho
 \\[1mm]
 &= \frac{2^m}{(2\pi)^{N}} \int_{\mathbb{R}^{N}} |\hat{\eta_\delta}(\rho)|^{2} \Big[(\log |\xi|)^m + \sum \limits^{m}_{j=1}
 \binom{m}{j} (\log|\xi|)^{m-j} (\log( \frac{|\rho-\xi|}{|\xi|}))^{j} \Big] d\rho,
  \end{align*}
  where   $ \rho=\gamma+\xi$.

   Observe that
 $$\log(\frac{|\rho-\xi|}{|\xi|}) \leq \log(\frac{|\rho|+|\xi|}{|\xi|}) =\log(1+\frac{|\rho|}{|\xi|}) \leq \frac{|\rho|}{|\xi|},\quad  0<\frac{r}{1+r}<1 \ \ {\rm for}\ r>0$$
then if choose $M=N+m+1$, we have that
\begin{align*}
 (\mathcal{L}_{m}e_{\xi,y},e_{\xi,y})_{L^{2}(\Omega)}
 &\leq 2^m (\log |\xi|)^m + 2^m  T_N c_M \sum \limits^{m}_{j=1} \binom{m}{j} (\log|\xi|)^{m-j} |\xi|^{-j} \int^{+\infty}_{0} \frac{r^j}{(1+r)^M}  r^{N-1} dr\\
 &\leq 2^m (\log |\xi|)^m + 2^m  T_N c_M \sum \limits^{m}_{j=1} \binom{m}{j} (\log|\xi|)^{m-j} |\xi|^{-j} \int^{+\infty}_{0} \frac{1}{(1+r)^{m-j+2}} dr \\
 &\leq 2^m (\log |\xi|)^m + 2^m    T_N c_{M}   \sum \limits^{m}_{j=1} \frac{\binom{m}{j}}{m-j+1} (\log|\xi|)^{m-j} |\xi|^{-j} .
  \end{align*}
 Therefore we find
\begin{align*}
\sum_{k}(\lambda-\lambda_{m,k}(\Omega))_{+} &\geq \frac{1}{(2\pi)^{N}} |\Omega_{\delta}|
 \int_{\mathbb{R}^{N}}\Big(\lambda - 2^m (\log |\xi|)^m - 2^m   T_N c_{M} \sum \limits^{m}_{j=1}    \frac{\binom{m}{j}}{m-j+1} (\log|\xi|)^{m-j} |\xi|^{-j}\Big)_{+} d\xi
\\[1mm]& =  T_N |\Omega_{\delta}|     \int^{\infty}_{0}  \big( f_\lambda(r) \big)_{+}  r^{N-1} dr,
   \end{align*}
where we used the definition of $f_\lambda$ by (\ref{fun b1}) with
 $$ T_N=\frac{ \omega_{_N}}{(2\pi)^{N}}\,\quad{and}\quad \,  \beta_{m,j} =  T_N   c_{M} \frac{\binom{m}{j}}{m-j+1}. $$

 Recall that  $r_0=e^{\frac12 \lambda^{\frac1m}}>1$, then we have that
 \begin{align*}
 \sum_{k}(\lambda-\lambda_{m,k}(\Omega))_{+}
 &\geq   T_N |\Omega_{\delta}|     \int^{\infty}_{0}  \big( f_\lambda(r) \big)_{+}  r^{N-1} dr
  \\[1mm]&\geq   T_N |\Omega_{\delta}|     \int^{r_0}_1  \big(\lambda - 2^m (\log r)^m - 2^m   \sum \limits^{m}_{j=1}  \beta_{m,j} (\log r)^{m-j} r^{-j} \big)_{+}  r^{N-1} dr
 \\[1mm]
 &\geq   T_N |\Omega_{\delta}|     \int^{r_0}_1 \big( \lambda - 2^m (\log r)^m -D_m \big( (\log r)^{m-1} r^{-1} +r^{-1}\big)  \big) \,  r^{N-1} dr\\[1mm]
 &=  T_N |\Omega_{\delta}|   \bigg(  \frac1N  \lambda r^{N}\bigg|^{r_0}_1
 - \frac1N  2^m r^{N} \sum \limits^{m}_{j=0} (-1)^{j} \frac{A_{m,j}}{2^j N^{j}} (\log r)^{m-j} \bigg|^{r_0}_1\\
 & \qquad \qquad \qquad \  -  D_m \frac{1}{N-1} \Big(r^{N-1} +r^{N-1}  \sum^{m-1}_{l=0} (-1)^{l} \frac{A_{m-1,l}}{2^l N^{l}}
 (\log r)^{{m-l-1}} \Big) \bigg|^{r_0}_1 \bigg),
 \end{align*}
 where
 $$D_m=\Big(\sup_{2\leq j\leq m} \sup_{r>1}\frac{(\log r)^{m-j} r^{-j}}{r^{-1}+(\log r)^{m-1} r^{-1}}\Big) \big(2^m    \sum \limits^{m}_{j=1}  \beta_{m,j}\big).$$
 Observe that
   \begin{align*}
  \lambda r^{N}\bigg|^{r_0}_1
 -2^m r^{N} \sum \limits^{m}_{j=0} (-1)^{j} \frac{A_{m,j}}{2^j N^{j}} (\log r)^{m-j} \bigg|^{r_0}_1
 = r_0^{N} \sum^{m}_{j=1}  (-1)^{j+1} \frac{A_{m,j}}{N^{j}}  \lambda ^{\frac{m-j}{m}} - \lambda
  + 2^m    (-1)^{m} \frac{A_{m,m}}{2^m N^{m}}
  \end{align*}
 and
\begin{align*}
& \quad   \Big(r^{N-1}+  r^{N-1}  \sum^{m-1}_{l=0} (-1)^{l} \frac{A_{m-1,l}}{2^l N^{l}}
 (\log r)^{{m-l-1}} \Big)  \bigg|^ {r_0}_1
\\[2mm] &=  r_0^{N-1} -1+   r_0^{N-1}  \sum^{m-1}_{l=0} (-1)^{l} \frac{A_{m-1,l}}{2^l N^{l}}
 (\log r_0)^{{m-l-1}}  -     (-1)^{m-1} \frac{A_{m-1,m-1}}{2^{m-1} N^{m-1}}
.
 \end{align*}
Note $$\lim_{\lambda \to +\infty} \frac{r_0^{N} \sum^{m}_{j=1}  (-1)^{j+1} \frac{A_{m,j}}{N^{j}}  \lambda ^{\frac{m-j}{m}} - \lambda
  + 2^m    (-1)^{m} \frac{A_{m,m}}{2^m N^{m}} }
{ r_0^{N}  \frac{A_{m,1}}{N}  \lambda ^{\frac{m-1}{m}} } = 1  $$
and
$$\lim_{\lambda \to +\infty} \frac{  r_0^{N-1} -1+ r_0^{N-1}  \sum^{m-1}_{l=0} (-1)^{l} \frac{A_{m-1,l}}{2^l N^{l}}
 (\log r_0)^{{m-l-1}}  -     (-1)^{m-1} \frac{A_{m-1,m-1}}{2^{m-1} N^{m-1}}}
{r_0^{N} \frac{A_{m,1}}{N}  \lambda ^{\frac{m-1}{m}} } =  0. $$

As a consequence, we obtain that
$$
 \liminf_{\lambda \to +\infty}  \frac{\sum_{k}(\lambda-\lambda_{m,k}(\Omega))_{+}}{\frac{A_{m,1}}{N} e^{\frac N2 \lambda^{\frac1m}}    \lambda ^{\frac{m-1}{m}}}
 \geq \frac { T_N |\Omega_{\delta}|}{N} , $$
which implies (\ref{es 2.2}), thanks to the arbitrariness of $\delta > 0$. \hfill$\Box$\medskip

\noindent {\bf Proof of Theorem \ref{teo 1}. } It follows by (\ref{re2.5}) and Lemma \ref{pr 2.2}   that
 \begin{align}\label{limt 4.1}
  \lim_{\lambda\to\infty} \Big(  \frac{2m}{N^2} \lambda^{\frac{m-1}{m}}\Big)^{-1}
  e^{-\frac N2 \lambda^{\frac{1}{m}}} \sum_{k}(\lambda-\lambda_{m,k}(\Omega))_{+} =  T_N |\Omega|.
  \end{align}


For $h>0$, we see that
 \begin{align}
  \frac{(\lambda+h-\lambda_{m,k}(\Omega))_{+} - (\lambda-\lambda_{m,k}(\Omega))_{+}}{h} \geq 1_{(-\infty,\lambda)}(\lambda_{m,k}(\Omega))  \label{bbb2.11}
  \end{align}
and
 \begin{align}
  \frac{(\lambda-\lambda_{m,k}(\Omega))_{+} - (\lambda-h-\lambda_{m,k}(\Omega))_{+}}{h} \leq 1_{(-\infty,\lambda)}(\lambda_{m,k}(\Omega)).\label{bbb2.12}
  \end{align}
From (\ref{limt 4.1}), for given $h>0$, there exists $\lambda^{*}>0$ such that  for $\lambda>\lambda^{*}$,
 \begin{align}\label{iq2.12}
 \Big(  \frac{2m}{N^2} \lambda^{\frac{m-1}{m}}\Big)^{-1} e^{-\frac N2 \lambda^{\frac{1}{m}}}\sum_{k}(\lambda-\lambda_{m,k}(\Omega))_{+}  \leq   T_N |\Omega|
 +\frac {h^2}{2},
 \end{align}

 \begin{align}\label{iq2.13}
 \Big(  \frac{2m}{N^2} \lambda^{\frac{m-1}{m}}\Big)^{-1} e^{-\frac N2 \lambda^{\frac{1}{m}}}\sum_{k}(\lambda-\lambda_{m,k}(\Omega))_{+}  \geq   T_N |\Omega|
 -\frac {h^2}{2} .
 \end{align}

 By (\ref{bbb2.11}), (\ref{iq2.12}) and (\ref{iq2.13}), we deduce that    for every  $h>0$
 $$\cN_{m}(\lambda) \leq \sum_{k} \frac{(\lambda+h-\lambda_{m,k}(\Omega))_{+} - (\lambda-\lambda_{m,k}(\Omega))_{+}}{h},$$
 for simplicity, we use the notations $\lambda_{m,k}=\lambda_{m,k}(\Omega)$,
 then for $\lambda>\lambda^{*}$ and $h>0$ small enough, we have
 \begin{align*}
 e^{-\frac N2 \lambda^{\frac{1}{m}}} \cN_{m}(\lambda)
 &\leq \frac{2m}{N^2} \lambda^{\frac{m-1}{m}}\,  \frac{\sum_{k} \frac{(\lambda+h-\lambda_{m,k})_{+} -(\lambda-\lambda_{m,k})_{+}}{h}}
 {\Big(  \frac{2m}{N^2} \lambda^{\frac{m-1}{m}}\Big)  e^{\frac N2 \lambda^{\frac{1}{m}}}}\\[1mm]
 &= \frac{1}{h}   \frac{2m}{N^2} \lambda^{\frac{m-1}{m}} \,
  \bigg[ \frac{\Big(  \frac{2m}{N^2} (\lambda+h)^{\frac{m-1}{m}}\Big) e^{\frac N2 (\lambda+h)^{\frac{1}{m}}}}{ \Big(  \frac{2m}{N^2} \lambda^{\frac{m-1}{m}}\Big) e^{\frac N2 \lambda^{\frac{1}{m}}} } \frac{\sum_{k} \Big((\lambda+h-\lambda_{m,k})_{+}}{\Big( \frac{2m}{N^2} (\lambda+h)^{\frac{m-1}{m}}\Big)  e^{\frac N2 (\lambda+h)^{\frac{1}{m}}}} -  \frac{\sum_{k}(\lambda-\lambda_{m,k})_{+}}{\Big(  \frac{2m}{N^2} \lambda^{\frac{m-1}{m}}\Big) e^{\frac N2 \lambda^{\frac{1}{m}}}} \bigg] \\[1mm]
  &\leq \frac{1}{h}\Big(  \frac{2m}{N^2} \lambda^{\frac{m-1}{m}}\Big)
  \Big[ \frac{\Big(  \frac{2m}{N^2} (\lambda+h)^{\frac{m-1}{m}}\Big) e^{\frac N2 (\lambda+h)^{\frac{1}{m}}}}{ \Big(  \frac{2m}{N^2} \lambda^{\frac{m-1}{m}}\Big) e^{\frac N2 \lambda^{\frac{1}{m}}} }  ( T_N |\Omega|+\frac {h^2}{2}) -  ( T_N |\Omega| - \frac {h^2}{2}) \Big].
 \end{align*}
For $\lambda>1$, denote
$$\Psi_\lambda (h)=e^{\frac N2\big( (\lambda+h)^{\frac{1}{m}}- \lambda^{\frac{1}{m}}\big)}
  (1+\frac{h}{\lambda})^{\frac{m-1}{m}}, $$
then $h\in(-1,1) \mapsto \Psi_\lambda (h)$ is continuous, differentiable,
 $\Psi_\lambda (0)=1$  for any $\lambda>0$   and
 $$
  \Psi_\lambda' (h)=\frac 1m e^{\frac N2 \big( (\lambda+h)^{\frac{1}{m}}- \lambda^{\frac{1}{m}}\big)}
  \Big(\frac N2  \lambda^{\frac{1-m}{m}} + (m-1) \lambda^{-1} (1+\frac{h}{\lambda})^{-\frac{1}{m}}\Big)
 $$
 and
 \begin{align*}
  \Psi_\lambda' (0)=  \frac 1m
  \Big(\frac N2  \lambda^{\frac{1-m}{m}} + (m-1) \lambda^{-1} \Big).
 \end{align*}
 Now we obtain that for any $h\in(0,1)$,
 \begin{align*}
  e^{-\frac N2 \lambda^{\frac{1}{m}}} \cN_{m}(\lambda)
 &\leq    \Big( \frac{2m}{N^2}  \lambda^{\frac{m-1}{m}}\, \frac{   \Psi_\lambda (h)  ( T_N |\Omega|+\frac {h^2}{2}) -  ( T_N |\Omega| - \frac {h^2}{2}) }
 {h}  \Big) \\[1mm]
 &=   \frac{2m}{N^2}  \lambda^{\frac{m-1}{m}}\,  \Big( T_N |\Omega|    \frac{  \Psi_\lambda (h)   -\Psi_\lambda (0)    }
 {h}+\Psi_\lambda (h) \frac h2    +\frac h2\Big).
 \end{align*}
 By the arbitrary of $h>0$, we can see that
 \begin{align*}
 \limsup_{\lambda\to\infty} e^{-\frac N2 \lambda^{\frac{1}{m}}} \cN_{m}(\lambda)
 &\leq    \limsup_{\lambda \to \infty}  \lim_{h\to0^+} \frac{2m}{N^2}  \lambda^{\frac{m-1}{m}}\,  \Big( T_N |\Omega|    \frac{  \Psi_\lambda (h)   -\Psi_\lambda (0)    }
 {h}+\Psi_\lambda (h) \frac h2    +\frac h2\Big)
 \\[1mm]
 &\leq  T_N |\Omega|  \limsup_{\lambda \to \infty} \Big( \frac{2m}{N^2}  \lambda^{\frac{m-1}{m}}  \Psi_\lambda' (0) \Big)\\
 &= \frac{  T_N |\Omega|}{N}.
  \end{align*}
That is
 \begin{align}\label{upp 1}
 \limsup_{\lambda\to\infty} e^{-\frac N2 \lambda^{\frac{1}{2}}} \cN_{m}(\lambda)
  \leq  \frac{  T_N |\Omega|}{N}  .
  \end{align}

 By (\ref{bbb2.12}), (\ref{iq2.12}) and (\ref{iq2.13}), we obtain that for every  $h>0$
 \begin{align*}
 e^{-\frac N2 \lambda^{\frac{1}{m}}} \cN_{m}(\lambda)
 &\geq   \frac{2m}{N^2} \lambda^{\frac{m-1}{m}}\,  \frac{\sum_{k} \frac{(\lambda-\lambda_{m,k})_{+} -(\lambda-h-\lambda_{m,k})_{+}}{h}}
 {\Big(  \frac{2m}{N^2} \lambda^{\frac{m-1}{m}}\Big)  e^{\frac N2 \lambda^{\frac{1}{m}}}}
 \\[1mm] &= \frac{1}{h}  \frac{2m}{N^2} \lambda^{\frac{m-1}{m}}\,
 \Big(\frac{\sum_{k}(\lambda-\lambda_{m,k})_{+}}{  \frac{2m}{N^2} \lambda^{\frac{m-1}{m}} e^{\frac N2 \lambda^{\frac{1}{m}}}}
 -  \frac{  \frac{2m}{N^2} (\lambda-h)^{\frac{m-1}{m}}  e^{\frac N2 (\lambda-h)^{\frac{1}{m}}}}{    \frac{2m}{N^2} \lambda^{\frac{m-1}{m}} e^{\frac N2 \lambda^{\frac{1}{m}}} } \frac{\sum_{k} (\lambda-h-\lambda_{m,k})_{+}}{  \frac{2m}{N^2} (\lambda-h)^{\frac{m-1}{m}}   e^{\frac N2 (\lambda-h)^{\frac{1}{m}}}}\Big)
 \\[1mm] &\geq\frac{1}{h}  \frac{2m}{N^2} \lambda^{\frac{m-1}{m}}\,
 \Big( T_N |\Omega| - \frac {h^2}{2}
 -  \frac{\Big(  \frac{2m}{N^2} (\lambda-h)^{\frac{m-1}{m}}\Big) e^{\frac N2 (\lambda-h)^{\frac{1}{m}}}}{ \Big(  \frac{2m}{N^2} \lambda^{\frac{m-1}{m}}\Big) e^{\frac N2 \lambda^{\frac{1}{m}}} }   ( T_N |\Omega|+\frac {h^2}{2}) \Big) ,
 \end{align*}
 then passing to the limit, we obtain that
 \begin{align*}
 \liminf_{\lambda\to+\infty} e^{-\frac N2 \lambda^{\frac{1}{m}}} \cN_{m}(\lambda)
 &\geq  \liminf_{\lambda\to+\infty}  \lim_{h\to0^+} \Big( \frac{2m}{N^2}  \lambda^{\frac{m-1}{m}}  \frac{    T_N |\Omega| - \frac {h^2}{2} - e^{\frac N2\big( (\lambda-h)^{\frac{1}{m}}- \lambda^{\frac{1}{m}}\big)}
  (1-\frac{h}{\lambda})^{\frac{m-1}{m}} ( T_N |\Omega|+\frac {h^2}{2})  }
 {h} \Big)  \\
 &= \liminf_{\lambda \to +\infty}  \lim_{h\to0^+}  \Big(  \frac{2m}{N^2} \lambda^{\frac{m-1}{m}}\Big)  \Big( T_N |\Omega| \frac{ \Psi_\lambda (0)   -\Psi_\lambda (-h)   }
 {h } -\frac h2 - \Psi_\lambda (-h) \frac h2\Big)
 \\[1mm]
 &\geq   \liminf_{\lambda \to +\infty}   \Big(   T_N |\Omega|  \frac{2m}{N^2} \lambda^{\frac{m-1}{m}}\Big) \Psi_\lambda' (0) \\[1mm]
 &=\frac{  T_N |\Omega|}{N}.
  \end{align*}
That is
 \begin{align}\label{upp 2}
 \liminf_{\lambda\to\infty} e^{-\frac N2  \lambda^{\frac{1}{m}}} \cN_{m}(\lambda)
 \geq \frac{  T_N |\Omega|}{N}.
 \end{align}
As a consequence, (\ref{limt 4.0}) follows by (\ref{upp 1}) and (\ref{upp 2}).  \hfill$\Box$\medskip

 \setcounter{equation}{0}
 \section{The lower bound of $\lambda_{m,k}(\Omega)$}
 \noindent

 \begin{lemma}\label{lm3.1}
 Let $m \geq 2$, and $a\geq \frac{2(m-1)}{N}$, then
 $$0\leq  \sum^{m}_{j=1} (-1)^{j+1} \frac{A_{m,j}}{N^{j+1}} a^{m-j} \leq  \frac{A_{m,1}}{N^{2}}a^{m-1} .$$
 \end{lemma}
 {\bf Proof. }Let $b_j= \frac{A_{m,j}}{N^{j+1}} a^{m-j} >0$, then
 $$\frac{b_j}{b_{j+1}} =  \frac{aN}{2(m-j)} \geq \frac{m-1}{m-j} \geq 1 ,\ \ {\rm for}\ j=1,\ 2,...,\ m-1. $$
Then for even $m$
 \begin{align*}
  0 \leq \sum^{m}_{j=1} (-1)^{j+1} b_j&=(b_1-b_2)+...+(b_{m-1}-b_m)
 \\& =b_1+(-b_2+b_3)+...+(-b_{m-2}+b_{m-1})-b_m \leq b_1\end{align*}
 and
 for odd   $m$,
 $$ 0\leq \sum^{m}_{j=1} (-1)^{j+1} b_j=b_1+(-b_2+b_3)+...+(-b_{m-1}+b_{m}) \leq b_1.$$
 We complete the proof.\hfill$\Box$\medskip

 \begin{lemma}\label{lm3.2}
 For  $m \geq 2$, we set
 $$h_m(\tau)= (\frac{\tau}{3})^{\frac 1{m-1}} - \log(\tau + e)\quad {\rm for} \ \, \tau > 0, $$
 then there exists $\tilde \tau_m > 1$, such that
 $$h_m(\tilde \tau_m)=0  \quad  {\rm and}\quad  h_m(\tau)>0  \ \ \ {\rm for}\   \tau> \tilde \tau_m .$$
  \end{lemma}
  {\bf Proof. } Note that
  $$\lim_{\tau \to \infty} \frac{\log(\tau + e)}{(\frac{\tau}{3})^{\frac 1{m-1}}} = \lim_{\tau \to \infty} \frac{\frac{1}{\tau + e}}{\frac 1{3m-3}  (\frac{\tau}{3})^{\frac{1}{m-1}-1}}
  =\lim_{\tau \to \infty} \frac{\frac{\tau}{3}}{\tau+e}  \frac{3m-3}{(\frac{\tau}{3})^{\frac{1}{m-1}}}   = 0,$$
  that is,   $$h_m(\tau) \to +\infty \ \ \ {\rm as}\ \tau\to +\infty.$$
   Moreover,
    $$h_m(1)= (\frac 13)^{\frac 1{m-1}} -  \log(1+e)<0. $$
  Thus there exists $\tilde  \tau_m > 1$, such that
 $$h_m(\tilde \tau_m)=0, \ {\rm and}\ h_m(\tau)>0  \ \ \ {\rm for}\   \tau >\tilde  \tau_m.$$
We complete the proof.\hfill$\Box$\medskip

For  $m \geq 2$, let
 $$f_1(t)=  (\frac N2)^{m-1}  e^{\frac{N}{2}t}  t^{m-1} \quad {\rm for}\ \, t\geq0, $$
which is continuous and  strictly increasing in $(0,+\infty)$ and $f_1(0)=0$. Then  for $\tau >0$,  there exists unique zero point $t_\tau$ of $f_1$, i.e.
$f_1(t_\tau)=\tau$.

  \begin{lemma}\label{lm3.3}
For $\tau >0$, let $t_\tau$ be the unique zero point  of $f_1$,
 then
 \begin{align}\label{aux 3.1}
 t_\tau\geq  \frac 2N \min\Big\{ \big(\frac{\tau}{e^{\tau_0}} \big)^{\frac 1{m-1}},\ \log \frac{\tau + \tau_m}{2(\log (\tau+e))^{m-1}}  \Big\},
   \end{align}
   where
   $$\tau_m=\max\{e,\tilde \tau_m\}\geq e, $$
   $\tilde \tau_m$ is from Lemma \ref{lm3.2} and $\tau_0\geq 1$ such that
   $$  \tau_0 ^{m-1}  e^{\tau_0}= \tau_m. $$
 \end{lemma}
 {\bf Proof. } Since $f_1$ is continuous and increasing  in $  (0,+\infty)$, the lower bound of $t_\tau$ could be done by constructing suitable $\bar t_\tau>0$ such that
 $$f_1(\bar t_\tau)\leq \tau.$$
From  Lemma \ref{lm3.2},  for $\tau>\tau_m  > 1$, one has
 $$3<   \frac{\tau}{(\log(\tau+e))^{m-1}} \leq  \tau,$$
for $\tau+\tau_m > \tau_m > 1$, we also have
$$\ 3<   \frac{\tau+\tau_m}{(\log(\tau+\tau_m+e))^{m-1}} \leq  \tau+\tau_m.$$

For  $\tau_0\geq 1$, let
$$\bar t_\tau= \left\{
  \begin{array}{lll}
\frac 2N \big(\frac{\tau}{e^{\tau_0}} \big)^{\frac 1{m-1}} \quad  &{\rm for}\ \, \tau\in(0,\tau_m],
 \\[2.5mm] \phantom{   }
\frac 2N  \log \frac{\tau}{(\log (\tau+e))^{m-1}} \quad    &{\rm for}  \ \, \tau>\tau_m,

  \end{array}
 \right.
 $$

 For $\tau > \tau_m$,

 \begin{align*}
 f_1(\bar t_\tau )
 &=\frac{\tau}{(\log(\tau+e))^{m-1}}   \Big(\log \frac{\tau}{(\log(\tau+e))^{m-1}}   \Big)^{m-1}  \\
 &\leq \tau \Big(\frac{\log \tau}{\log(\tau+e)}    \Big)^{m-1} \\
 &< \tau.
  \end{align*}

 For $ 0<\tau  \leq \tau_m  =  \tau_0 ^{m-1}  e^{\tau_0}$, we see that
 $$(\frac{\tau}{e^{\tau_0}})^{\frac 1{m-1}}  -  \tau_0 \leq 0,$$
 then
 $$f_1 (\bar t_\tau)= \tau  e^{(\frac{\tau}{e^{\tau_0}})^{\frac 1{m-1}}}    e^{-\tau_0}  =  \tau  e^{(\frac{\tau}{e^{\tau_0}})^{\frac 1{m-1}}  -  \tau_0}    \leq   \tau,$$

 Thus, we have that
  \begin{align}\label{aux 3.2}
  \bar t_\tau\geq \left\{
  \begin{array}{lll}
\frac 2N \big(\frac{\tau}{e^{\tau_0}} \big)^{\frac 1{m-1}} \quad  &{\rm for}\ \, \tau\in(0,\tau_m],
 \\[3mm] \phantom{   }
\frac 2N  \log \frac{\tau}{(\log (\tau+e))^{m-1}} \quad    &{\rm for}  \ \, \tau>\tau_m.
  \end{array}
 \right.
  \end{align}
 For $\tau> \tau_m$, we have that
 $$\frac{\tau}{(\log (\tau+e))^{m-1}} > \frac{\tau + \tau_m}{2(\log (\tau+e))^{m-1}} >\frac{\tau + \tau_m}{2(\log (\tau+ \tau_m +e))^{m-1}} > \frac 32,$$
 thus
 $$\log \frac{\tau}{(\log (\tau+e))^{m-1}} > \log \frac{\tau + \tau_m}{2(\log (\tau+e))^{m-1}}>0.$$
 Moreover, we see that for $\tau\in[0, \tau_m]$,
  \begin{align}\label{aux 112}
  \frac{\tau + \tau_m}{2(\log (\tau+e))^{m-1}}>\frac{\tau + \tau_m}{2(\log (\tau+\tau_m+e))^{m-1}}\geq\frac32>1,    
  \end{align}
 we can extend  (\ref{aux 3.2}) to
  \begin{align}\label{aux 3.3}
 \bar t_\tau \geq  \frac 2N \min\Big\{ \big(\frac{\tau}{e^{\tau_0}} \big)^{\frac 1{m-1}},\ \log \frac{\tau + \tau_m}{2(\log (\tau+e))^{m-1}}  \Big\}  .
  \end{align}
Thus the Lemma \ref{lm3.3} is proved.\hfill$\Box$\medskip

\noindent  {\bf Proof of Theorem \ref{teo 2}. }
It infers by Proposition \ref{pr 2.1} that for  $m \geq 2$,
  \begin{align}\label{aux 3.4}
  \cN_{m}(\lambda) \leq   T_N |\Omega| \frac{1}{\eta-\lambda} e^{\frac N2\eta^{\frac{1}{m}}}
 (\sum \limits^{m}_{j=1} (-1)^{j+1} \frac{A_{m,j}}{N^{j+1}} \eta^{\frac{m-j}{m}})
 \end{align}
and we shall take the precise parameters $\lambda = \lambda_{m,k}(\Omega) + \frac1{N^m}$  and $\eta  = \lambda_{m,k}(\Omega) + P_m(\Omega)$,
where
 $$P_m(\Omega)=  \big(\lambda_{m,1}(\Omega)\big)_- +\frac{2^m}{N^m}  (m-1)^m + \frac1{N^m},$$
 here $$a_\pm=\max\{0,\pm a\}. $$
Note that $\lambda_{m,k}(\Omega) +\big(\lambda_{m,1}(\Omega)\big)_-\geq \lambda_{m,k}(\Omega) - \min\{0,\lambda_{m,1}(\Omega)\} \geq \lambda_{m,k}(\Omega) -  \lambda_{m,1}(\Omega) \geq0$, then
 $$\eta   =\lambda_{m,k}(\Omega) +\big(\lambda_{m,1}(\Omega)\big)_-+\frac{2^m}{N^m}  (m-1)^m +\frac1{N^m}
 > \max\{\lambda,\ \frac{2^m}{N^m}  (m-1)^m   \} > 0.$$
 For simplicity, we use the notations
 $$\lambda_{m,k}=\lambda_{m,k}(\Omega),\quad \lambda_{m,1}=\lambda_{m,1}(\Omega),\quad P_m=P_m(\Omega).$$
 then
 $$\frac{A_{m,j}}{N^{j+1}} (\lambda_{m,k}+P_m)^{\frac{m-j}{m}} \geq
 \frac{A_{m,j+1}}{N^{j+2}} (\lambda_{m,k}+P_m)^{\frac{m-j-1}{m}}\ \ {\rm for}\ j=1,\ 2,\ ...,\ m-1, $$
and now we apply the Lemma \ref{lm3.1} to obtain that
 $$  \sum \limits^{m}_{j=1} (-1)^{j+1} \frac{A_{m,j}}{N^{j+1}} (\lambda_{m,k}+P_m)^{\frac{m-j}{m}} \leq \frac{A_{m,1}}{N^{2}} (\lambda_{m,k}+P_m)^{\frac{m-1}{m}}.$$
Thus, (\ref{aux 3.4}) becomes
 \begin{align*}
 k &\leq \frac{ T_N |\Omega|}{P_m} e^{\frac N2 (\lambda_{m,k}+P_m)^{\frac{1}{m}}}
 \Big(\sum \limits^{m}_{j=1} (-1)^{j+1} \frac{A_{m,j}}{N^{j+1}} (\lambda_{m,k}+P_m)^{\frac{m-j}{m}}\Big) \\
 &\leq \frac{ T_N |\Omega|}{P_m} e^{\frac N2 (\lambda_{m,k}+P_m)^{\frac{1}{m}}}  \frac{A_{m,1}}{N^2} (\lambda_{m,k}+P_m)^{\frac{m-1}{m}},
  \end{align*}
  that is,
  $$\frac{N^2 P_m k}{2m  T_N |\Omega|} \leq   e^{\frac N2 (\lambda_{m,k}+P_m)^{\frac{1}{m}}} (\lambda_{m,k}+P_m)^{\frac{m-1}{m}}$$
 or
  $$\frac{N^{m+1} P_m k}{2^m  m  T_N |\Omega|} \leq  (\frac N2)^{m-1} e^{\frac N2 (\lambda_{m,k}+P_m)^{\frac{1}{m}}} (\lambda_{m,k}+P_m)^{\frac{m-1}{m}}.$$

Recall that
  $$ \quad a_m= \frac{N^{m+1}}{ 2^m  m  } \frac{(2\pi)^{N}}{ \omega_{_N}},\quad b_m=  \frac{2^m}{N^m}  (m-1)^m +\frac1{N^m},\quad c_m=a_mb_m
  ,  \quad \tau_m = \tau_0 ^{m-1}  e^{\tau_0}.  $$
we apply  Lemma \ref{lm3.3} with
$$\tau = \frac{N^{m+1} P k}{2^m  m  T_N |\Omega|}=\frac{ a_m P_m k}{ |\Omega|}, \qquad  t = (\lambda_{m,k}+P_m)^{\frac{1}{m}}$$
 to obtain that
 \begin{align}\label{e 4.1}
(\lambda_{m,k}+P_m)^{\frac{1}{m}}  \geq \bar t_\tau \geq  \frac 2N \min\Big\{ \big(\frac{\tau}{e^{\tau_0}} \big)^{\frac 1{m-1}},\ \log \frac{\tau + \tau_m}{2(\log (\tau+e))^{m-1}}  \Big\}.
  \end{align}

When either  $m$ is odd and $\lambda_{m,1}(\Omega)\geq0$,  or $m$ is even,   we have that $\lambda_{m,1}(\Omega)\geq 0$, then $P_m=\frac{2^m}{N^m}  (m-1)^m +\frac1N$,
and (\ref{e 4.1}) implies
 \begin{align}
\lambda_{m,k}(\Omega)  \geq  &\big(\frac{2}{N}\big)^m   \min\Big\{  \big(  \frac{c_m}{e^{\tau_0}} \big)^{\frac m{m-1}}      ( \frac{ k }{|\Omega|})^{\frac m{m-1}}, \nonumber\\[2mm]
 &\qquad\qquad\quad\ \Big(\log \big(c_m  \frac{k }{|\Omega|}    + \tau_0 ^{m-1}  e^{\tau_0}\big)  - (m-1)\log \big(\log (c_m  \frac{k }{|\Omega|} + e)  \big) -\log2 \Big)^m    \Big\}   - b_m.\label{low 1.1}
  \end{align}
It infers by \cite[Theorem 1.3$(v)$]{C0} that  $\lambda_{m,1} \geq0$   when $\Omega\subset B_{r_0}$ for some $r_0>0$ small. \smallskip

When   $m$ is odd and $\lambda_{m,1}<0$,  then  $P_m(\Omega)=-\lambda_{m,1}(\Omega)  +\frac{2^m}{N^m}  (m-1)^m +\frac1{N^m}$ and
    \begin{align}
\lambda_{m,k}  \geq  \big(\frac{2}{N}\big)^m   \min\Big\{ &\big(  \frac{a_m}{e^{\tau_0}} \big)^{\frac m{m-1}}   \big(\frac{ P_m  }{  |\Omega|} \big)^{\frac m{m-1}}   k^{\frac m{m-1}},\nonumber \\
 &\Big(\log (a_m \frac{P_m  }{|\Omega|}  k + \tau_0 ^{m-1}  e^{\tau_0})  - (m-1)\log  \big(\log (a_m  \frac{P_m }{|\Omega|}  k + e)  \big) - \log2 \Big)^m    \Big\}   - P_m(\Omega)\label{xxx-1}
  \end{align}
as claimed. \hfill$\Box$\medskip

 \begin{corollary}\label{cr 1.1}
 $(i)$ When  either  $m$ is odd and $\Omega\subset B_{r_0}(0)$ for some $r_0>0$ or
  $m$ is even,  then
  $$ \lambda_{m,1}(\Omega)\geq  0  $$
  and
 \begin{align}\label{cr low-1}
   \liminf_{|\Omega| \to0^+}\lambda_{m,1}(\Omega)  \Big(  \log \frac{1}{|\Omega|(-\log |\Omega| )   }  \Big)^{-m}\geq (\frac 2N)^{m}.
     \end{align}

  $(ii)$  When   $m$ is odd and $\lambda_{m,1}(\Omega)<0$,  then
 \begin{align}\label{cr low-2}
  \lambda_{m,1}(\Omega) \geq b_m- d_m  |\Omega|
   \end{align} \medskip
  where $d_m \geq \max \big\{\frac{e^{\tau_0}}{ a_m} \big( (m-1)^m + \frac1{2^m} \big)^{\frac{m-1}{m}}, \,
 \frac{1}{ a_m} \big( 2 e^{(m-1)^m + \frac1{2^m} } -  \tau_0 ^{m-1} e^{\tau_0}  \big) \big\} > 0. $

  \end{corollary}

\noindent{\bf Proof. }
 Recall that
  $$ \quad a_m= \frac{N^{m+1}}{ 2^m  m  } \frac{(2\pi)^{N}}{ \omega_{_N}},\quad b_m=  \frac{2^m}{N^m}  (m-1)^m +\frac1{N^m},\quad c_m=a_m b_m.  $$
 {\it Part 1:}   When  either  $m$ is odd and $\Omega\subset B_{r_0}(0)$ for some $r_0>0$ or
  $m$ is even,   we have $\lambda_{m,1}(\Omega)\geq 0$.   From the inequality (\ref{low 1.1}) with $k=1$,  we have that
  \begin{align*}
   \lambda_{m,1}(\Omega)  > \big(\frac{2}{N}\big)^m \min \bigg\{\big(  \frac{c_m}{e^{\tau_0}} \big)^{\frac m{m-1}}   \big(\frac{ 1  }{  |\Omega|} \big)^{\frac m{m-1}}   ,
 \Big(\log \frac{  \frac{c_m  }{|\Omega|}  + \tau_0 ^{m-1}  e^{\tau_0}  }{ 2 \big(\log (\frac{c_m  }{|\Omega|}  + e) \big)^{m-1}}  \Big)^m  \bigg\}
    - b_m.
    \end{align*}
As $|\Omega|\to0^+$,
  \begin{align}\label{low-aab}
   \lambda_{m,1}(\Omega)  &> \big(\frac{2}{N}\big)^m  \Big(\log \frac{  \frac{c_m  }{|\Omega|}  + \tau_0 ^{m-1}  e^{\tau_0}  }{ 2 \big(\log (\frac{c_m  }{|\Omega|}  + e) \big)^{m-1}} \Big)^m -b_m,
    \end{align}
  which implies (\ref{cr low-1}). \smallskip

 {\it Part 2:}  we deal the case that   $m$ is odd and $\lambda_{m,1}(\Omega) <0$.   Then $P_m(\Omega)=-\lambda_{m,1}(\Omega)  +b_m $
 and (\ref{xxx-1}) implies that
     \begin{align*}
 b_m\geq  \big(\frac{2}{N}\big)^m   \min\Big\{ &\big(  \frac{a_m}{e^{\tau_0}} \big)^{\frac m{m-1}}   \big(\frac{ P_m  }{  |\Omega|} \big)^{\frac m{m-1}}   ,\nonumber \\
 &\Big(\log (a_m \frac{P_m  }{|\Omega|}   + \tau_0 ^{m-1}  e^{\tau_0})  - (m-1)\log \big(\log (a_m  \frac{P_m }{|\Omega|}   + e)  \big) - \log2 \Big)^m    \Big\},
  \end{align*}
 then
 $$b_m\geq \frac{2^m}{N^m} \big(  \frac{a_m}{e^{\tau_0}} \big)^{\frac m{m-1}}   \big(\frac{ P_m  }{  |\Omega|} \big)^{\frac m{m-1}}   $$
 or
 $$
 b_m
 \geq \frac{2^m}{N^m}  \Big(\log \frac{  a_m \frac{P_m  }{|\Omega|}  + \tau_0 ^{m-1}  e^{\tau_0}  }{ 2}  \Big)^m
 \geq  \frac{2^m}{N^m} \Big(\log \frac{  a_m \frac{P_m  }{|\Omega|}  + \tau_0 ^{m-1}  e^{\tau_0}  }{ 2\big(\log (a_m  \frac{P_m  }{|\Omega|} + e) \big)^{m-1}}  \Big)^m ,
 $$
 thus
 $$P_m(\Omega) \leq \big( \frac N2  \big)^{m-1}  \frac{e^{\tau_0}}{ a_m} b_m^{\frac{m-1}{m}}  |\Omega|$$
 or
 $$P_m (\Omega)\leq \frac{1}{ a_m} \big(2 e^{\frac N2  b_m^{\frac1m} } -  \tau_0 ^{m-1} e^{\tau_0}  \big)  |\Omega|
 \leq \frac{1}{ a_m} \big( 2 e^{\frac N2  b_m^{\frac1m} } -  e  \big)  |\Omega|,$$
 then there is
 $$
 d_m \geq \max \big\{\big( \frac N2  \big)^{m-1}  \frac{e^{\tau_0}}{ a_m} b_m^{\frac{m-1}{m}},\
 \frac{1}{ a_m} \big(2 e^{\frac N2  b_m^{\frac1m} } -  e  \big)   \big\} > 0,
 $$
 there holds
$$
 \lambda_{m,1}(\Omega) \geq b_m - d_m  |\Omega|.
 $$
   We complete the proof. \hfill$\Box$\medskip

\begin{remark}
When $m$ is even and $|\Omega|$ large, 
 we see that 
 $$\big( \frac{c_m}{e^{\tau_0}} \big)^{\frac m{m-1}}   \big(\frac{ 1  }{  |\Omega|} \big)^{\frac m{m-1}}\to 0  \ \ \ {\rm as}\  |\Omega| \to \infty$$
 and
 $$\log \frac{  \frac{c_m  }{|\Omega|}  + \tau_0 ^{m-1}  e^{\tau_0}  }{ 2 \big(\log (\frac{c_m  }{|\Omega|}  + e) \big)^{m-1}} \geq
 \log  \frac32>0 $$
 by (\ref{aux 112}). 
 
Thus, we see that $|\Omega|$ large
$$ \lambda_{m,1}(\Omega) > - b_m, $$
which is useless, because we have $\lambda_{m,1}(\Omega)\geq0$. 
\end{remark}

  \begin{remark}
 When $m=1$,  from  \cite[Corollary 1.4]{CV0}
  \begin{equation}\label{Weyl-upplow}
  \frac{2}{N} \Big(\ln k+\ln \frac{2}{eN d_N|\Omega|}\Big)\leq  \lambda_{1,k}(\Omega)\leq  \frac{2}{N} \ln k+c_0\big(\ln\ln (k+e)\big)^2+\frac2N\ln \frac{|B_1|}{|\Omega|},
    \end{equation}
  where $d_N>0$ is some constant depending on dimension $N$.
  Particular for   $k=1$,
   \begin{equation}\label{Weyl-upplow}
  \frac{2}{N}  \ln \frac{2}{eN d_N}-\ln |\Omega| \leq  \lambda_{1,1}(\Omega)\leq   c_0\big(\ln\ln (1+e)\big)^2+\frac2N\ln  |B_1| -\ln |\Omega|.
    \end{equation}
 Compared with (\ref{Weyl-upplow}),     the upper bound   (\ref{cr low-1}) is sharp, but
the lower bound (\ref{cr low-2}) is still rough.
 \end{remark}
    \bigskip\bigskip

   \noindent{\bf \small Acknowledgements:} {\footnotesize H. Chen  is supported by NSFC,  No. 12071189, 12361043,  by Jiangxi Province Science Fund No. 20232ACB201001. L. Chen is supported by Postgraduate Innovation Fund of Educational Department in Jiangxi Normal University YJS2023075.}


\medskip

\begin{thebibliography}{99}
\bibitem{BCH} M. Bhakta, H. Chen, H. Hajaiej: On the bounds of the sum of eigenvalues for a Dirichlet problem involving mixed fractional Laplacians.  {\it J. Diff.  Eq. 317},    1--31 (2022).


 \bibitem {CS0}   L. Caffarelli,   S. Salsa,    L. Silvestre:
Regularity estimates for the solution and the free boundary to the
obstacle problem for the fractional Laplacian. {\it Invent.
Math. 171,} 425--461 (2008).



\bibitem{CS1}  L. Caffarelli,   L.  Silvestre:  An extension problem related to the fractional Laplacian.
{\it Comm. Part. Diff. Eq. 32}, 1245--1260 (2007).



\bibitem {CL} H. Chen, P. Luo:  Lower bounds of Dirichlet eigenvalues for some degenerate elliptic operators.  {\it  Calc.Var. Part. Diff. Eq. 54(3),} 2831--2852 (2015).



 \bibitem {CW0}  H.Y.  Chen,   T. Weth: The Dirichlet problem for the logarithmic Laplacian.
{\it Comm. Part. Diff. Eq. 44(11)},   1100--1139 (2019).

 \bibitem {CV0} H.Y. Chen,  L. V\'eron:  Bounds for eigenvalues of the Dirichlet problem for the logarithmic Laplacian. {\it Adv.    Calc.  Var. 16(3)},     541--558  (2023).

 \bibitem {C0} H.Y. Chen:  Taylor's expansions of Riesz convolution and the fractional Laplacians with respect to the order.  {\it arXiv: 2307.06198v2 }

 \bibitem {CS} Z. Chen, R. Song:  Two-sided eigenvalue estimates for subordinate
processes in domains. {\it J. Funct. Anal. 226},    90--113 (2005).

 \bibitem {CgY} Q.  Cheng, H. Yang:  Bounds on eigenvalues of Dirichlet Laplacian.  {\it   Math. Ann. 337,} 159--175   (2007).


\bibitem {CgW} Q. Cheng, G. Wei:  A lower bound for eigenvalues of a clamped plate problem.  {\it  Calc.Var. Part. Diff. Eq. 42(3/4),} 579--590 (2011).

\bibitem {CDK} V.  Crismale, L. De Luca, A. Kubin, M. Ponsiglione:
The variational approach to  $s$-fractional heat flows and the limit cases  $s\to0^+$  and  $s\to 1^-$.
 {\it J. Funct. Anal. 284(8)},  Paper No. 109851, 38 pp.  (2023).

 \bibitem {EGE}  E. Di Nezza,  G. Palatucci,   E.  Valdinoci:  Hitchhiker's guide to the fractional Sobolev spaces.
{\it  Bull.   Sci. Math.   136(5)},  521--573 (2012).

\bibitem {FJW} P.  Feulefack, S. Sven and T. Weth:
Small order asymptotics of the Dirichlet eigenvalue problem for the fractional Laplacian.
{\it J. Fourier Anal. Appl. 28(2)}, No. 18, 44 pp. (2022).

\bibitem {Ge} L. Geisinger:
A short proof of Weyl's law for fractional differential operators. {\it J. Math. Phys. 55(1)},  011504, 7 pp. (2014).

 \bibitem {HY} E. M. Harrell II, S. Y. Yolcu:  Eigenvalue inequalities for Klein-Gordon operators.
 {\it  J. Funct. Anal. 256(12),} 3977--3995 (2009).

\bibitem {HS} V. Hern\'andez Santamar\'ia,  A. Saldana:  Small order asymptotics for nonlinear fractional problems.
{\it Calc. Var. PDE 61(3)},   Paper No. 92, 26 pp. (2022).


 \bibitem {KM}  M. Kassmann, A.  Mimica: Intrinsic scaling properties for nonlocal operators. {\it J. Eur. Math. Soc. 19(4)}, 983--1011 (2013).

\bibitem {K} P. Kr\"oger:  Estimates for sums of eigenvalues of the Laplacian.  {\it J. Funct. Anal. 126(1),} 217--227 (1994).


 \bibitem {LW}   A. Laptev,   T. Weth:  Spectral properties of the logarithmic Laplacian.
{\it Anal. Math. Phys. 11(3)},   no. 133, 24 pp. (2021).

\bibitem {LY} P. Li and  S.-T.Yau:  On the Schr\"odinger equation and the eigenvalue problem. {\it Commun. Math. Phys. 88(3),} 309--318 (1983).

  \bibitem {L} E. Lieb:   The number of bound states of one-body Schr\"odinger operators and the Weyl problem.  {\it Proc. Sym. Pure Math. 36,} 241--252 (1980).



    \bibitem {M}  A. Melas:  A lower bound for sums of eigenvalues of the Laplacian. {\it Proc. Am. Math. Soc. 131(2),} 631--636 (2003).



\bibitem{MN}  R. Musina,    A. Nazarov:  On fractional Laplacians. {\it Comm. Part. Diff. Eq.  39,} 1780--1790  (2014).


   \bibitem {P} G. P\'olya:   On the Eigenvalues of Vibrating Membranes
 (In Memoriam Hermann Weyl).  {\it  Proc. Lond. Math. Soc. 3(1),} 419--433 (1961).

 \bibitem{Ro} G. Rozenbljum: Distribution of the discrete spectrum of singular operator. {\it Dokl. Akad. Nauk
SSSR 202}, 1012--1015 (1972).




 \bibitem {RS1}  X. Ros-Oton,   J. Serra:  The Pohozaev identity for the fractional Laplacian.  {\it Arch.   Ration. Mech.   Anal. 213,} 587--628 (2014).



\bibitem {RS}  X. Ros-Oton,   J. Serra:  The Dirichlet problem for the fractional
laplacian: regularity up to the boundary. {\it J. Math. Pures Appl. 101},  275--302 (2014).




\bibitem {S07}  L. Silvestre:   Regularity of the obstacle problem for a fractional power of the Laplace operator.  {\it Comm. Pure   Appl. Math. 60 (1),} 67--112 (2007).



\bibitem {W} H. Weyl:  Das asymptotische Verteilungsgesetz der Eigenwerte linearer partieller Differentialgleichungen (mit einer Anwendung auf die Theorie der Hohlraumstrahlung). {\it Math. Ann. 71(4),} 441--479 (1912).


\bibitem {YY} S. Yolcu, T. Yolcu:   Estimates for the sums of eigenvalues of the fractional Laplacian on a bounded domain. {\it Commun. Contemp. Math. 15(3),} 1250048 (2013).


\end{thebibliography}
\end{document}